\documentclass[]{amsart}

 \newtheorem{thm}{Theorem}[section]
 \newtheorem{cor}[thm]{Corollary}
 \newtheorem{lem}[thm]{Lemma}
 \newtheorem{prop}[thm]{Proposition}
 \theoremstyle{definition}
 \newtheorem{defn}[thm]{Definition}
 \newtheorem{defns}[thm]{Definitions}
 \theoremstyle{remark}
 \newtheorem{rem}[thm]{Remark}
 \newtheorem{rems}[thm]{Remarks}
 \newtheorem{exam}[thm]{Example}
 
 \newtheorem{conj}[thm]{Conjecture}
 \newtheorem{que}[thm]{Question}
 \newtheorem{rem and def}[thm]{Remark and Definition}
 \newtheorem{def and rem}[thm]{Definition and Remark}
 
 \numberwithin{equation}{section}
 \usepackage [all,2cell]{xy}

\begin{document}

\title{ Cluster Structure on Generalized Weyl Algebras }
\date{}

\author{Ibrahim Saleh}
% -----------------------------------------------------------
\maketitle
% -----------------------------------------------------------

\begin{abstract} We introduce  a  class of non-commutative algebras that carry  non-commutative   cluster structure which are generated by identical copies of  generalized Weyl algebras.  Equivalent conditions for the finiteness of the set of the cluster variables of these cluster structures are provided.  Mutations along with some  combinatorial data, called \textit{cluster strands,} arising from the cluster structure are used to construct representations of generalized Weyl algebras.

\end{abstract}
% -----------------------------------------------------------
\maketitle
% -----------------------------------------------------------
\tableofcontents

\section{Introduction} Cluster algebras were introduced by S. Fomin and A. Zelevinsky in [8, 9, 2, 10, 17]. A cluster algebra is a commutative algebra with a distinguished set of generators called cluster variables and particular type of relations called mutations. A quantum  version was introduced in  [3] and [5, 6, 7]. The original motivation was to create a combinatorial algebraic framework to study total positivity and dual canonical basis in coordinate rings of certain semisimple algebraic groups.

Generalized Weyl algebras  were first introduced by  V. Bavula in [1] and separately  as  Hyperbolic algebras by  A. Rosenberg in [14]. Their  motivation was to find a ring theoretical frame work to study the representation theory of some important ``small algebras"  such as the first Heisenberg algebra, Weyl algebras, the universal enveloping algebra of the Lie algebra $sl(2)$. A complete list of ``small algebras" can be found in [14]. Also, in [14] Rosenberg has obtained the representation theory of all ``small algebras" using the Hyperbolic algebra as a frame work.

 In this paper, we show that by relaxing the commutativity between cluster variables and some frozen variables (coefficients variables) we can extend the theory of cluster algebras to include some non-commutative algebras which are generated by isomorphic copies of  generalized Weyl algebras. To achieve this goal, we introduced particular non-commutative seed-like combinatorial data called \emph{presseds}, each preseed of rank $n$ is defined by iteration from a rank one  preseed. Every Fomin-Zelevinsky (coefficient free) rank one seed $(\{x\}, \cdot_{x})$ gives rise to a preseed  of rank one by attaching a valued star quiver with center at the  vertex $\cdot_{x}$ and assigning a set of frozen variables, one frozen variable at each vertex of the star quiver. Here the frozen variables associated with the exchange variable $x$ do not necessary commute with it. Every  preseed of rank $n$ is defined through  an increasing set of $n-1$ (nested)  preseeds  of ranks $1,\ldots, n-1$ respectively, Definition 3.2.  A valued star quiver is called \emph{balanced} if  the (componentwise) sum of the valuations of the arrows point in toward the center vertex equals the sum of the valuations  of the arrows  point out equals $(a, a)$, for some non-negative integer $a$. A preseed is called balanced if each of its star quivers is balanced.

  The set of all cluster variables produced from a preseed  is not necessarily finite, even if the underlying quiver is of Dynkin type, Examples 3.9 and 3.10. In this paper we provide equivalent  conditions on a preseed  for its  set of cluster variables to be finite, such as in Theorems 3.15 and Corollary 4.9 which are rephrased respectively as follows

 \begin{thm} Let $p$ be a balanced preseed  in the ambient division ring $\mathcal{D}$. If $\phi$ is a $\mathcal{D}$-automorphism that fixes the frozen variables  such that   for every frozen variable $f$  associated to the initial cluster  variable $x$ we have
  \begin{equation}\label{}
   \nonumber  fx=\phi (x)f.
  \end{equation}
 Then the set of all cluster variables of $p$ is finite if and only if $\phi$ is of finite order.
 \end{thm}

 \begin{cor} Let $p$ be a balanced preseed  with a non zero element $q$ in the field $K$ such that for every initial cluster  variable $x$, we have
  \begin{equation}\label{}
   \nonumber  fx=qxf, \text{for each frozen variable $f$  associated to}\ x.
  \end{equation}
  Then the set of all cluster variables of $p$ is finite if and only if $q$ is an $m^{th}$-root of unity, for some natural number $m$.
  \end{cor}
 Although Corollary 1.2 can be seen as a consequence of Theorem 1.1, in this paper  we provide independent proofs for both of them. \\
Every generalized Weyl algebra of rank $n$ gives rise to a preseed  $p$ of rank $n$ endowed with an automorphism $\theta$ over the coefficients ring (the group ring of the group generated by the frozen variables), Example 4.5. In  Theorem 4.12 we show that the associated  cluster algebra $\mathcal{H}(p)$ satisfies the following

\begin{enumerate}
  \item  The algebra $\mathcal{H}(p)$  is generated by (possibly) infinite isomorphic copies of the associated generalized Weyl algebra, each vertex in the exchange graph of $p_{n}$ gives rise to two copies of them;
  \item There are  $n$ rank one preseeds $p_{1}(x_{1}), \ldots, p_{1}(x_{n})$ such that
  \begin{equation}\label{}
\mathcal{H}(p)=\mathcal{H}(p_{1}(x_{1}))\otimes\cdots\otimes\mathcal{H}(p_{1}(x_{n})).
\end{equation}

\end{enumerate}

Let $V_{n}$ be the $K$-left span of the cluster monomials of $\mathcal{H}(p)$. In  Definition 5.5 we use right and left mutations, given in Definition 3.3, to  introduce an action of generalized Weyl algebra in $V_{n}$. The combinatorial structure of the cluster monomials gives rise to combinatorial datum called \emph{cluster strands}, which are particular elements of  $V_{n}$,  Definition 5.8. Some properties of the cluster strands are provided in Lemma 5.12. The submodules generated by cluster strands are called \emph{strand submodules}. The properties of the strand submodules are studied in Proposition 5.15 and 5.16 and Corollary 5.17.

 \begin{conj} Strand submodules are indecomposable.
 \end{conj}

The paper is organized as follows. Section $2$ is devoted to basic definitions of cluster algebras associated with valued quivers. In Section 3, we introduce the notion of preseeds   and their mutations. Examples and properties of preseeds  are also given. In the same section we provide equivalent conditions for a preseed  to produce a finite set of cluster variables, Theorems 3.15. In Theorem 3.17, we introduce a class of $\mathcal{D}$-automorphisms that preserve the set of cluster variables. Weyl cluster algebras are defined in Section 4. The main results of Section 4 are Corollary 4.9 and Theorem 4.12 which give some basic properties of Weyl cluster algebras. Section 5 is where we introduce an action of generalized Weyl algebras on the space of cluster monomials. In the same section we introduce the  cluster strands. Some of their basic properties are in Lemma 5.12. Some  Properties of strand submodules are given in Proposition 5.15 and 5.16.

 Through out the paper, $K$ is a field of zero characteristic and the notation $[1, k]$ stands for the set $\{1,\ldots ,k\}$.

%%% ----------------------------------------------------------------------

\section{Cluster algebras associated with valued quivers}
For more details about the material of this section  refer to  [16, 17, 11, 2, 8].
\subsection{Valued quivers}
 \begin{itemize}
\item
   \emph{A valued quiver} of rank $n$ is a quadruple $Q=(Q_{0}, Q_{1}, V, d)$, where
  \begin{itemize}
 \item $Q_{0}$ is a  set of $n$ vertices labeled by numbers from the set $[1, n]$;
  \item  $Q_{1}$ is called the \emph{set of arrows} of $Q$ and consists of  ordered pairs of vertices, that is $Q_{1}\subset Q_{0}\times Q_{0}$;
  \item $V$ is a function $V:Q_{1}\rightarrow \mathbb{N}\times\mathbb{N}$, $(i, j)\mapsto(v_{ij},v_{ji})$, $V$ is called the \emph{valuation} of $Q$;
  \item $d=(d_{1},\cdots, d_{n})$, where $d_{i}$ is a positive integer for each $i$, such that $d_{i}v_{ij}=v_{ji}d_{j}$, for every $i, j\in Q_{0}$.
   \end{itemize}
   In the case of $(i,j)\in Q_{1}$, then there is an arrow oriented from $i$ to $j$ and in notation we shall use the symbol $\xymatrix{{\cdot}_{i} \ar[r]^{(v_{ij},v_{ji})}&{\cdot}_{j}}$. If $v_{ij}=v_{ji}=1$ we simply write $\xymatrix{{\cdot}_{i} \ar[r]&{\cdot}_{j}}$.\\
In this paper, we moreover assume that $(i, i)\notin Q_{1}$ for every $i\in Q_{0}$, and  if $(i, j)\in Q_{1}$ then $(j, i)\notin Q_{1}$. The vector of positive integers $d=(d_{1},\cdots, d_{n})$ does not play any role in the context of this paper, so it will be ignored from now on.
 \item   If $v_{ij}=v_{ji}$ for every $(v_{ij},v_{ji})\in V$ then $\Gamma$ is called  \emph{equally  valued quiver}.
\item  We say that the valued quiver $\Gamma=(Q_{0}, Q_{1}, V)$ is \emph{connected}, if for every $v, v' \in {Q_{0}}$, there is a sequence of vertices $v=v_{1}, \cdots, v_{l}=v'$  such that for $t=1, \cdots, l-1$, either $(v_{t}, v_{t+1})$ or $(v_{t+1}, v_{t})$ is in $Q_{1}$, in other words,  any pair of subsequent vertices $v_{t}$ and $v_{t+1}$ are connected by an arrow.
 \end{itemize}

 \begin{rems}
 \begin{enumerate}
   \item
   Every (non valued) quiver $Q$ without loops nor $2$-cycles corresponds to an equally valued quiver which has an arrow $(i, j)$ if there is at least one arrow directed from $i$ to $j$ in $Q$ and with the valuation   $(v_{ij}, v_{ji})=(m, m)$, where $m$ is the number of arrows from $i$ to $j$.

  \item
Every  valued quiver of rank $n$ corresponds to a skew symmetrizable integer matrix   $B(Q)=(b_{ij})_{i,j\in[1,n]}$ given by
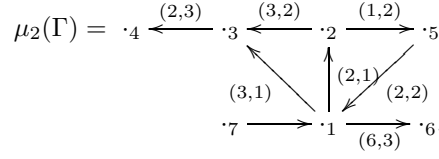
\begin{equation}\label{}
  b_{ij}=\begin{cases} v_{ij}, & \text{if} \ (i,j)\in Q_{1},\\
    0, & \text{if \  neither} \  (i,j)  \ \text{nor} \  (j,i) \  \text{is in} \  Q_{1},\\
-v_{ij}, & \text{ if }(j,i)\in Q_{1}.
    \end{cases}
\end{equation}
Conversely,  given a skew symmetrizable $n\times n$ matrix $B$, a valued quiver $Q_{B}$ can be easily defined  such that $B(Q_{B})=B$. This gives rise to a  bijection between the skew-symmetrizabke  $n\times n$ integral matrices $B$ and the valued quivers with set of vertices $[1,n]$, up to isomorphism fixing the vertices.

  \end{enumerate}
 \end{rems}

\begin{defn}[\emph{Valued quivers mutations}] Let $Q$ be a valued quiver. The mutation  $\mu_{k}(Q)$ at a vertex  $k$  is defined  through Fomin-Zelevinsky's mutation of the associated skew-symmetrizable matrix. The mutation of a skew symmetrizable matrix $B=(b_{ij})$ on the direction $k\in [1, n]$ is given by $\mu_{k}(B)=(b'_{ij})$, where
\begin{equation}\label{}
b'_{ij}=\begin{cases} -b_{ij}, & \text{if} \ k \in \{i,j\},\\
    b_{ij}+\text{sign}(b_{ik})\max(0, b_{ik}b_{kj}), & \text{otherwise.}
    \end{cases}
  \end{equation}

\end{defn}

\begin{rems}
\begin{enumerate}
  \item  Let  $Q=(Q_{0}, Q_{1}, V)$ be a valued quiver. The new valued quiver $\mu_{k}(Q)=(Q_{0}, Q'_{1}, V')$, obtained from $Q$ by applying mutation at the  vertex  $k$, can be described using the mutation of $B(Q)$ as follows:
       We obtain $Q'_{1}$ and $V'$  by altering  $Q_{1}$ and $V$, based on the following rules
\begin{enumerate}
 \item replace the  pairs $(i, k)$ and $(k,j)$  with $(k,i)$ and $(j,k)$  respectively and switch the components of the ordered pairs of their valuations;
  \item if  $(i,k), (k,j)\in Q_{1}$, such that at least $i$ or $j$ is in $Q_{0}$ but $(j,i)\notin Q_{1}$  and $(i,j)\notin Q_{1}$ (respectively $(i,j)\in Q_{1}$) add the  pair $(i, j)$ to $Q'_{1}$, and give it the valuation $(v_{ik}v_{kj},v_{ki}v_{jk})$ (respectively change its valuation to $(v_{ij}+v_{ik}v_{kj},v_{ji}+v_{ki}v_{jk})$);
 \item if $(i,k)$, $(k,j)$ and $(j, i)$ in $Q_{1}$, then we have three cases
 \begin{enumerate}
   \item if $v_{ik}v_{kj}<v_{ij}$, then keep $(j,i)$ and change its valuation to $(v_{ji}-v_{jk}v_{ki}, |-v_{ij}+v_{ik}v_{kj}|)$;
   \item if $v_{ik}v_{kj}>v_{ij}$, then replace $(j,i)$ with $(i,j)$ and change its valuation to $( -v_{ij}+v_{ik}v_{kj}, |v_{ji}-v_{jk}v_{ki}|)$;
   \item if $v_{ik}v_{kj}=v_{ij}$,  then remove $(j,i)$ and its valuation.
 \end{enumerate}
\end{enumerate}
 \item One can see that; $\mu^{2}_{k}(Q)=Q$ and  $\mu_{k}(B(Q))=B(\mu_{k}(Q))$ at each vertex $k$.

 \end{enumerate}
\end{rems}

     \begin{exam}Let
     \begin{equation}\label{}
     \Gamma=\xymatrix{
     \cdot_{4}& \cdot_{3}   \ar[l]_{(2,3)}\ar[r]^{(2,3)}&   \cdot_{2}\ar[d]^{(1,2)}&\cdot_{5}\ar[l]_{(2,1)}\\
  &\cdot_{7} \ar[r]&\cdot_{1}\ar[ul]^{(9,3)}\ar[r]^{(6,3)}&\cdot_{6}.}
     \end{equation}
     One can see that $\Gamma$ is a valued quiver with $d=(1, 2, 3, 2, 1, 2, 1 )$. Applying mutation at the vertex  $2$, produces the following valued quiver

  \begin{equation}\label{}
      \nonumber  \mu_{2}(\Gamma)=\xymatrix{      \cdot_{4}& \cdot_{3}   \ar[l]_{(2,3)}&   \cdot_{2} \ar[l]_{(3,2)}\ar[r]^{(1,2)}&\cdot_{5}\ar[dl]^{(2, 2)}\\
  &\cdot_{7} \ar[r]&\cdot_{1}\ar[ul]^{(3,1)}\ar[u]_{(2,1)}\ar[r]_{(6,3)}&\cdot_{6}.}
     \end{equation}

  \end{exam}

\subsection{Cluster algebras} [17] Let $\mathcal{F}$  be an ambient field of rational functions in $n$ independent  variables over $\mathbb{Q}(t_{1}, \ldots, t_{m})$. A \emph{seed}  in $\mathcal{F}$  is a pair $(X, Q)$, where
\begin{itemize}
  \item $X=\{x_{1},\ldots, x_{n}\}$ forms a free generating set of $\mathcal{F}$,  and
  \item $Q=(Q_{0}, Q_{1}, V)$ is a valued quiver with  $Q_{0}=\{1,\ldots, n, n+1,\ldots, n+m\}$, where  vertices $1,\ldots, n$ are called \emph{exchange vertices} and $n+1,\ldots, n+m$ are the called \emph{frozen vertices}.
\end{itemize}
The variables $x_{1},\ldots, x_{n}$  are associated with the exchange vertices and they are called \emph{exchange cluster variables} and  the variables $t_{1}, \ldots, t_{m}$ are associated with the frozen vertices and they are called \emph{frozen variables}.

\begin{defn}[Seed mutations] Let $p=(X, Q)$ be a seed in $\mathcal{F}$ and let $k\in [1,n]$. Applying the \emph{seed mutation} $\mu_{k}$  on $(X, Q)$ produces  a new seed $\mu_{k}(X, Q)=(\mu _{k}(X), \mu _{k}(Q))$ by  setting $\mu _{k}(X)=\{x_{1}, \ldots,x'_{k},\ldots, x_{n},  t_{n+1},\ldots, t_{n+m}\}$  where $x'_{k}$ is defined by the so-called \emph{exchange relations}:

\begin{equation}\label{}
   x'_{k}x_{k}=\prod_{j,\xymatrix{{\cdot}_{n+j} \longrightarrow {\cdot}_{k}}} t_{n+j}^{v_{n+j,k}}\prod_{i, \xymatrix{{\cdot}_{i} \longrightarrow{\cdot}_{k}}} x_{i}^{v_{ik}}+\prod_{j, \xymatrix{{\cdot}_{k}\longrightarrow{\cdot}_{n+j}}} t_{n+j}^{v_{k,n+j}}\prod_{i, \xymatrix{{\cdot}_{k} \longrightarrow{\cdot}_{i}}} x_{i}^{v_{ki}}.
\end{equation}
And $\mu_{k}(Q)$ is the mutation of $Q$ at the vertex $k$, given in Definition 2.3 and Remarks 2.4. The elements of $\mathcal{F}$ obtained by applying iterated mutations on the elements $\{x_{1}, \ldots, x_{n}\}$ are called \emph{cluster variables}.
   \end{defn}

\begin{defns}[Cluster algebra and exchange graph]\begin{enumerate}
\item
Let $\mathcal{X}$ be the set of all cluster variables of $\mathcal{F}$ produced from a seed $(X, Q)$. The  \emph{cluster algebra} $\mathcal{A}=\mathcal{A}(X, Q)$ is  the $\mathbb{Z}[\mathbb{P}]$-subalgebra of $\mathcal{F}$ generated by $\mathcal{X}$, where $\mathbb{P}$ is the (free) abelian group generated by the frozen variables written multiplicatively.

\item The \emph{exchange graph} of $\mathcal{A}(X, Q)$, denoted by $\mathbb{G}(X, Q)$, is the $n$-regular graph whose vertices are labeled by the seeds that can be obtained from $(X, Q)$ by applying some sequence of mutations, and whose edges correspond to mutations. Two adjacent seeds in $\mathbb{G}$ can be obtained from each other by applying a mutation $\mu_{k}$ for some $k\in [1, n]$.

\end{enumerate}
\end{defns}

\begin{thm} [8, Theorem 3.1, Laurent Phenomenon]  The cluster algebra
$\mathcal{A}(X, Q)$ is contained in the integral ring of Laurent polynomials $\mathbb{Z}[\mathbb{P}][x^{\pm}_{1},\ldots ,x^{\pm}_{n}]$.
 \end{thm}

\section{Preseeds}

Before introducing  \emph{preseeds}, we will  introduce an  increasing filtration of  division rings of fractions by iteration and a particular type of quivers known as \emph{star quivers}.\\ For each $m$ in $[1, n]$, let  $\mathbb{P}_{m}$ be a finitely  generated free abelian group, written multiplicatively,  with set of  generators
 \begin{equation}\label{}
   F^{m}=\bigcup^{m}_{i=1}F_{i} \ \ \text{where} \ \   F_{i}=\{f_{i1},\ldots,f_{im_{i}}\}.
 \end{equation}

 Let $\textit{R}_{1}=K[\mathbb{P}_{1}]$ be the  group ring of $\mathbb{P}_{1}$ over $K$. Let $D_{1}$ be an Ore domain containing $\textit{R}_{1}$ such that there is  $t_{1}\in D_{1}$ so that $\{t^{\alpha _{1}}_{1}; \alpha _{1}\in \mathbb{Z}\}$ form a basis for $D_{1}$ as a left $R_{1}$-module. Let $\mathcal{D}_{1}$ denote the set of right fractions $ab^{-1}$ with $a, b \in D_{1}$, and $b\neq 0$; two such fractions $ab^{-1}$ and $cd^{-1}$ are identified if $af = cg$ and $bf = dg$ for some non-zero $f, g \in D_{1}$. The ring $D_{1}$ is embedded into $\mathcal{D}_{1}$ via $d \mapsto d\cdot 1^{-1}$. The addition and multiplication in $D_{1}$ extend to $\mathcal{D}_{1}$ so that $\mathcal{D}_{1}$ becomes a division ring. (Indeed, we can define $ab^{-1}+ cd^{-1}=(ae + cf )g^{-1}$ where non-zero elements $e, f$ , and $g$ of $D_{1}$ are chosen so that $be = df = g$; similarly,
$ab^{-1} \cdot cd^{-1}=ae\cdot(df)^{-1}$, where non-zero $e, f \in D_{1}$ are chosen so that $cf = be$). In such case we say $\mathcal{D}_{1}$ is the division ring of fractions in $t_{1}$ of $D_{1}$ over $R_{1}$. Now, for $i\in [2, n]$, let $R_{i}=K[\mathbb{P}_{i}]$ and  $D_{i}$ be an Ore domain containing (as sub rings) $\textit{R}_{i}$ and $D_{i-1}$ such that there is  $t_{i}\in D_{i}$ so that
 \begin{equation}\label{}
   t_{i}t_{j}=t_{j}t_{i} \ \ \text{and} \ \ t_{i}f_{jr}=f_{jr}t_{i},\ \ \text{for every} \  i,j\in [1, n], j<i, \  \text{for all} \ \ r\in[1, m_{j}];
 \end{equation}
 and $\{t^{\alpha _{i}}_{i}; \alpha _{i}\in \mathbb{Z}\}$ form a basis for $D_{i}$ as a left $R_{i}$-module. Let $\mathcal{D}_{i}$ be  the division ring of fractions in $t_{i}$ of $D_{i}$ over $R_{i}$. For each $i\in [1, n]$, the elements of the set $F_{i}$ are called \emph{frozen variables}. More details about Ore domains can be found in [13] and [2]. The following graph is meant to help readers understand the relations between the rings $\textit{R}_{i}, D_{i}$ and $\mathcal{D}_{i}, i=1,\cdots, n$.

 \begin{equation}\label{}
  \nonumber \textit{R}_{1}\subset \textit{R}_{2} \subset \cdots \subset  \textit{R}_{n}\\
   \end{equation}
   \begin{equation}\label{}
  \nonumber  \cap \    \  \   \   \  \cap \ \ \ \ \  \cdots  \ \ \ \ \cap\\
    \end{equation}
    \begin{equation}\label{}
  \nonumber  D_{1} \subset D_{2} \subset \cdots \subset D_{n}\\
    \end{equation}
    \begin{equation}\label{}
  \nonumber  \cap \    \  \   \   \  \cap \ \ \ \ \  \cdots \ \ \ \ \cap\\
    \end{equation}
    \begin{equation}\label{}
  \nonumber \mathcal{D}_{1}\subset\mathcal{D}_{2} \subset \cdots \subset\mathcal{D}_{n}.
 \end{equation}

\

\begin{defns}[{Valued star quivers}]
\begin{itemize}
  \item A valued quiver  $\Gamma=(Q_{0}, Q_{1}, V)$  is called a \emph{valued star quiver} with center at $k\in Q_{0}$ if we have
      \begin{equation}\label{}
        \nonumber Q_{1}\subset (\{k\}\times Q_{0})\cup (Q_{0}\times \{k\}).
      \end{equation}
      Furthermore, $\Gamma$ is called a balanced star quiver if
      \begin{equation}\label{}
      (\sum_{j,\xymatrix{{\cdot}_{k} \longrightarrow {\cdot}_{j}}}v_{kj}, \sum_{j,\xymatrix{{\cdot}_{k} \longrightarrow {\cdot}_{j}}}v_{jk})=(\sum_{i,\xymatrix{{\cdot}_{i} \longrightarrow {\cdot}_{k}}}v_{ik}, \sum_{i,\xymatrix{{\cdot}_{i} \longrightarrow {\cdot}_{k}}}v_{ki})=(a_{k}, a_{k}),
      \end{equation}
      and in this case the non-negative integer   $a_{k}$ is called the \emph{frozen component} of $\Gamma$.
      \item A set of $n$ star quivers $\Gamma=\{\Gamma_{1}, \ldots, \Gamma_{n}\}$ is said be  \emph{balanced} with \emph{frozen rank}   $(a_{1}, \ldots, a_{n})$  if each  $\Gamma_{k}$  is balanced with frozen component $a_{k}, k=1,\ldots, n$.

      \end{itemize}
      \end{defns}
The following valued quiver is an example of a balanced valued star quiver of frozen rank $9$ and $d=(6, 4, 12, 24, 6, 3)$

  \begin{equation}\label{}
   \nonumber  \Gamma=\xymatrix{
     \cdot_{3}   \ar[r]^{(2,4)}&   \cdot_{1}\ar[dl]_{(4,1)}\ar[dr]^{(3, 6)}\ar[d]&\cdot_{2}\ar[l]_{(6,4)}\\
   \cdot_{4}&\cdot_{5} &\cdot_{6}.}
     \end{equation}

Although, in general, presseds can be defined using any valued quiver,  here they are  defined using valued star quivers which serves best the purpose of the paper. From now on we will omit the word valued from the term  valued star quiver.

\begin{defn}[{Preseeds}]  \begin{itemize}
\item
A \emph{preseed } $p_{1}$ of rank $1$ in  $\mathcal{D}_{1}$ is the  triple  $(\{F_{1}\}, \{x_{1}\}, \{\Gamma_{1}\})$, where
 \begin{enumerate}
 \item $F_{1}$ is as described in (3.1);

  \item $x_{1}$ is an element of  $\mathcal{D}_{1}$  such that there is an $\textit{R}_{1}$- linear automorphism on $\mathcal{D}_{1}$ that
 fixes the frozen variables and sends $t_{1}$ to $x_{1}$. The element  $x_{1}$ is  called an \emph{exchange cluster variable} and the set

 \begin{equation}\label{}
 \nonumber   \widetilde{X}:=\{f_{11},\ldots,f_{1m_{1}},x_{1}\}
 \end{equation}
 is called the \emph{extended cluster} of $p_{1}$;

 \item $\Gamma _{1}$ is a star quiver of rank $m_{1}+1$. The center vertex  $\cdot_{1}$ of $\Gamma _{1}$   is called \emph{exchange vertex}  and all other vertices are called \emph{frozen vertices}.

 \end{enumerate}

 \item A \emph{preseed } $p_{n}$ of rank $n$ in  $\mathcal{D}_{n}$ is the  triple $(F,  X, \Gamma)$, where $F=\{F_{1}, \ldots, F_{n}\}$ (as given in (3.1)),  $X=\{x_{1},\ldots, x_{n}\}$ and $\Gamma=\{\Gamma_{1}, \ldots, \Gamma_{n}\}$ such that $p_{k}=(\{ F_{k}\}, \{ x_{k}\}, \{ \Gamma_{k}\}\}$ is  a preseed  of rank $1$ in $\mathcal{D}_{k}$, for every $k \in [1, n]$.  The following set

      \begin{equation}\label{}
   \nonumber \widetilde{X}= \{f_{11},\ldots, f_{1m_{1}}, \ldots, f_{n1},\ldots, f_{nm_{n}}, x_{1}, \ldots,  x_{n}\}
 \end{equation}
is called the \emph{extended cluster} of $p_{n}$. Furthermore, $p_{n}$ is called \emph{balanced preseed}  if $\Gamma$ is a balanced set of star quivers and  the \emph{frozen rank} of $p_{n}$ is the same as the frozen rank  of $\Gamma$.

\end{itemize}

\end{defn}

\begin{defn}[Preseeds   mutations] Let $p_{n}=(F, X, \Gamma)$  be a preiseed    in  $\mathcal{D}_{n}$. For each  $k\in [1,n]$, two new triples  $\mu^{R}_{k}(p_{n})=(F, \mu^{R}_{k}(X), \mu_{k}(\Gamma))$ and  $\mu^{L}_{k}(p_{n})=(F, \mu^{L}_{k}(X), \mu_{k}\Gamma))$ can be obtained from $p_{n}$ as follows

\begin{itemize}
  \item (Right mutation)
   \begin{equation}\label{}
  \mu^{R}_{k}(x_{i})=\left\{
    \begin{array}{ll}
      \Large(\prod_{j, i_{j}\rightarrow i} f_{ij}^{v_{i_{j}i}}\prod_{j\rightarrow i} x_{j}^{v_{ji}}+\prod_{j, i\rightarrow i_{j}} f_{ij}^{v_{ii_{j}}}\prod_{j,i\rightarrow j} x_{j}^{v_{ij}}\Large)x^{-1}_{i}, & i=k; \\
      x_{i}, & i\neq k.
    \end{array}
  \right.
 \end{equation}

  \item (Left mutation)

   \begin{equation}\label{}
  \mu^{L}_{k}(x_{i})=\left\{
    \begin{array}{ll}
      x^{-1}_{i} \Large(\prod_{j, j\rightarrow i} x_{j}^{v_{ji}}\prod_{j, i_{j}\rightarrow i} f_{ij}^{v_{i_{j}i}}+\prod_{j, i\rightarrow j} x_{j}^{v_{ij}}\prod_{j, i\rightarrow i_{j}} f_{ij}^{v_{ii_{j}}}\Large), & i=k; \\
      x_{i}, & i\neq k.
    \end{array}
  \right.
 \end{equation}
    \item The   mutation $\mu_{k}(\Gamma)$  is as defined in Definition 2.2 and Remarks 2.3.

\end{itemize}

   \end{defn}

\begin{prop} Let $p_{n}=(F, X, \Gamma)$ be a  preseed    in  $\mathcal{D}_{n}$. Then the following are true
\begin{enumerate}
\item
  For any sequence of right mutations (respectively left)  $\mu^{R}_{i_{1}}\mu^{R}_{i_{2}}\ldots \mu^{R}_{i_{q}}$, we have
     $\mu^{R}_{i_{1}}\mu^{R}_{i_{2}}\ldots \mu^{R}_{i_{q}}(p_{n})$  (respectively $\mu^{L}_{i_{1}}\mu^{L}_{i_{2}}\ldots \mu^{L}_{i_{q}}(p_{n})$) is again a  preseed .
 \item For every $k \in [1, n]$,
 \begin{equation}\label{}
    \mu^{R}_{k} \mu^{L}_{k}(p_{n})=\mu^{L}_{k}\mu^{R}_{k}(p_{n})=p_{n}.
 \end{equation}
\end{enumerate}
\end{prop}

\begin{proof} We  prove part $(1)$ for  $\mu^{R}_{k}(p_{n})$, (respectively $\mu^{L}_{k}(p_{n})$) and the proof for an arbitrary sequence of right
(respectively left) mutations is by induction on the length of the sequence. From (3.4) (respectively (3.5)) one has $\mu^{R}_{k}(x_{k})$ (respectively $\mu^{L}_{k}(x_{k})$) is an expression in the elements of the set $\{ x^{-1}_{k}\} \cup F_{k}$. Then (3.2)  guarantees that $\{x_{1}, \ldots, x_{k-1}, \mu^{R}_{k}(x_{k}),  x_{k+1}, \ldots, x_{n}\}$ (respectively $\{x_{1},\ldots,x_{k-1}, \mu^{L}_{k}(x_{k}),  x_{k+1}, \ldots, x_{n}\}$) is a commutative set. The commutativity of  the elements of the set  $\{\mu^{R}_{k}(x_{k})\} \cup F_{j}$ (respectively the elements of the set $\mu^{L}_{k}(x_{k})\cup F_{j}$) for $j\neq k$ is again due to that the expression of  $\mu^{R}_{k}(x_{k})$  (respectively $\mu^{L}_{k}(x_{k})$) contains only elements of $\{ x^{-1}_{k}\} \cup F_{k}$ which is by (3.2) commute with elements of $F_{j}$.  Part $(2)$ is immediate using (3.2), (3.4) and (3.5) and the fact that mutation is involutive on valued quivers.
\end{proof}

\begin{defn}[Cluster sets] Let $p_{n}$ be a preseed  in $\mathcal{D}_{n}$. An element $y \in \mathcal{D}_{n}$ is said to be a \emph{ cluster variable}  if $y$  is  a cluster variable in some  seed $q_{n}$, where  $q_{n}$ is obtained from  $p_{n}$  by applying some sequence of (right or left) mutations. The set of all  cluster variables of $p_{n}$  is called   \textit{the  cluster set} of $p_{n}$  and is denoted by  $\mathcal{X }(p_{n})$.
\end{defn}

\begin{rem}
\begin{enumerate}
  \item  From the definition of preseeds,  each  exchange vertex  is connected only to its associated frozen vertices. Then from the proof of Part 1 of  Proposition 3.4, one can see that every cluster variable in $\mathcal{D}_{n}$, can be written as a  Laurent  expression  in exactly one  cluster variable  and the frozen variables associated to it in some pressed. Which is a major difference between cluster variables produced from preseeds and cluster variables produced from other  non-commutative seeds such as  quantum seeds introduced in [3].
  \item   Mutations of preseeds  are not involutive but they are invertible, in the  sense of Part 2 Proposition 3.4,   however  mutations of  classical or quantum  seeds are involutive.

\end{enumerate}

\end{rem}

\begin{defn}
 A quadruple $(F, X, \Gamma,  \varphi)$  is said to be $\varphi$-commutative preseed in $\mathcal{D}_{n}$ if $(F, X, \Gamma)$ is a preiseed    and $\varphi$ is  an $R_{n}$-linear automorphism   of $\mathcal{D}_{n}$, such that the following equations are satisfied

\begin{equation}
fx_{i}=\varphi (x_{i})f,  \ \  \forall f \in F_{i}, \ \  \forall i \in [1,n].
\end{equation}

\end{defn}
One can see equations (3.7) induces the equations
\begin{equation}
f^{a}x_{i}=\varphi^{a} (x_{i})f^{a}, \forall f \in F_{i}, \forall i \in [1,n], a \in \mathbb{Z}_{\geq 0}.
\end{equation}
And
   \begin{equation}
  x_{i}f^{a}=f^{a} \varphi^{-a} (x_{i}), \forall f \in F_{i}, \forall i \in [1,n], a \in \mathbb{Z}_{\geq 0}.
   \end{equation}
An example of a $\varphi$-commutative preseed is given in next section, Example 4.5.
\begin{prop} The  properties of balanced  and  $\varphi$-commutative of preseeds are invariant under preseeds   mutations.
\end{prop}
\begin{proof} One can see that the mutation of balanced preseed is again a balanced preseed with the same frozen rank.\\Now  we  show that $\varphi$-commutativity of preseeds is invariant under right mutations and for left mutation is quite similar. For every $k\in[1, n]$, the right mutation $\mu^{R}_{k}$  of $p_{n}$ gives rise to an $R_{n}$-automorphism  $\psi:\mathcal {D}_{n}\rightarrow \mathcal {D}_{n}$ induced by
 \begin{equation}\label{}
       \psi(x_{j}) :=\mu^{R}_{k}(x_{j}),  \ \forall j\in [1,n].
 \end{equation}
We will show that $\mu^{R}_{k}(p_{n})$ is $\psi \varphi \psi^{-1}$-commutative. Let $i\in[1, n]$, $f\in F_{i}$. We have

\begin{eqnarray}
% \nonumber to remove numbering (before each equation)
 \nonumber f\mu_{k}(x_{i}) &=& \psi(f x_{i}) \\
 \nonumber  &=&  \psi(\varphi (x_{i}f))\\
 \nonumber  &=&  \psi\varphi  \psi ^{-1}(\mu_{k}(x_{i}))f.
   \end{eqnarray}
\end{proof}
Let $p=(\{x_{1}\}, \cdot_{x_{1}})$ be a coefficient free seed of rank $1$ in the field of fractions $K(t)$. In this case there is only one more  seed   $(\{\frac{2}{x_{1}}\}, \cdot_{\frac{2}{x_{1}}})$ which is mutation equivalent to $p$.  The Fomin-Zelevinsky (commutative) cluster algebra of $p$ is the algebra of polynomials with integral coefficients $\mathcal{A}=\mathbb{Z}[x_{1}, \frac{2}{x_{1}}]$. In the following  we will see  two examples of attaching star quivers at vertex  $\cdot_{x_{1}}$ to produce preseeds.

\begin{exam} \textbf{The simplest non balanced   preseed }. Let $p_{1}$ be the  seed $(\{F_{1}\},\{x_{1}\},\{\Gamma_{1}\})$ where $F_{1}=\{f_{11}\}$ and $\Gamma_{1}$ is the following star quiver
 \begin{equation}\label{}
\nonumber     \xymatrix{
 \cdot_{f_{11}} &  \cdot_{x_{1}} \ar[l] }.
\end{equation}

Applying mutation at the vertex  $\cdot_{x_{1}}$, we obtain the following cluster variables

\begin{eqnarray*}
% \nonumber to remove numbering (before each equation)
 x_{1}&\stackrel{\mu^{L} _{k}}{\Rightarrow}& x_{1}^{-1}(f_{11}+1)\\
   &\stackrel{\mu^{L} _{k}}{\Rightarrow}& (f_{11}+1)^{-1}x_{1}(f_{11}+1) \\
      &\stackrel{\mu^{L} _{k}}{\Rightarrow}& (f_{11}+1)^{-1}x^{-1}_{1}(f_{11}+1)^{2} \\
      &\stackrel{\mu^{L} _{k}}{\Rightarrow}& (f_{11}+1)^{-2}x_{1}(f_{11}+1)^{+2} \\
   &\ldots &\\
  &\stackrel{\mu^{L} _{k}}{\Rightarrow}&(f_{11}+1)^{-k}x^{-1}_{1}(f_{11}+1)^{k+1}\\
  &\stackrel{\mu^{L} _{k}}{\Rightarrow}&(f_{11}+1)^{-(k+1)}x_{1}(f_{11}+1)^{k+1}\\
   &\ldots &,\\
\end{eqnarray*}
and
\begin{eqnarray*}
% \nonumber to remove numbering (before each equation)
   x_{1}&\stackrel{\mu^{R} _{k}}{\Rightarrow}& (f_{11}+1)x_{1}^{-1} \\
   &\stackrel{\mu^{R} _{k}}{\Rightarrow}& (f_{11}+1)x_{1}(f_{11}+1)^{-1} \\
      &\stackrel{\mu^{R} _{k}}{\Rightarrow}& (f_{11}+1)^{2}x^{-1}_{1}(f_{11}+1)^{-1} \\
      &\stackrel{\mu^{R} _{k}}{\Rightarrow}& (f_{11}+1)^{2}x_{1}(f_{11}+1)^{-2} \\   &\ldots &\\
  &\stackrel{\mu^{R} _{k}}{\Rightarrow}&(f_{11}+1)^{k+1}x^{-1}_{1}(f_{11}+1)^{-k}\\
  &\stackrel{\mu^{R} _{k}}{\Rightarrow}&(f_{11}+1)^{k+1}x_{1}(f_{11}+1)^{-(k+1)}\\
   &\ldots &.\\
\end{eqnarray*}

Then, we have the  infinite cluster  set  $\mathcal{X } (p_{1})=\{x_{1}, (1+f_{11})^{k+1}x_{1}^{-1}(1+f_{11})^{-k}, (1+f_{11})^{k}x_{1}(1+f_{11})^{-k},(1+f_{11})^{-k}x_{1}^{-1}(1+f_{11})^{k+1},(1+f_{11})^{-k}x_{1}^{-1}(1+f_{11})^{k}, k \in \mathbb{Z}\}$.
 Later in this article, we will see that this seed is related to first Weyl algebra.
\end{exam}

\begin{exam} \textbf{The cluster set of the simplest (nontrivial) balanced $\varphi$-commutative preseed}. Consider the  seed $p_{1}=(\{F_{1}\},\{x_{1}\},\{\Gamma_{1}\},\varphi)$
where $F_{1}=\{f_{11},f_{12}\}$ and $\Gamma_{1} $ is the following star quiver with frozen rank is (1)
\begin{equation}\label{}
\nonumber     \xymatrix{
 \cdot_{f_{11}} &  \cdot_{x_{1}} \ar[l]   &  \ar[l]\cdot_{f_{12}}}.
\end{equation}
If $\varphi $ be a  $\textit{R}$-linear automorphism of $\mathcal{D}_{1}$
satisfying the conditions (3.7). Then this seed produces the cluster set $\mathcal{X } (p_{1})$ given by
\begin{equation}\label{}
   \nonumber \{(f_{11}+f_{12})x_{1}^{-1}, x_{1}^{-1}(f_{11}+f_{12}), \varphi ^{k} (x_{1}),  (f_{11}+f_{12})\varphi^{-k} (x_{1}),  \varphi^{-k} (x_{1})(f_{11}+f_{12}); k\in \mathbb{Z} \}.
\end{equation}
One can see that $\mathcal{X } (p_{1})$ is a finite set if and only if $\varphi$
is of finite order.
\end{exam}

\begin{rem} Examples 3.9 and 3.10 show that the Fomin-Zelevinsky finite type classification [9] does not work in the preseed  case in general.
\end{rem}

\begin{lem} Let $p_{n}=(F, X, \Gamma, \varphi)$  be a  $\varphi$-commutative preseed. If $\Gamma_{k}$ is a balanced star quiver, then  we have

\begin{equation}\label{}
    (\mu^{R}_{ k})^{2}(x_{k})=\varphi ^{a_{k}} (x_{k}) \ \ \text{for some nonnegative integer} \ \ a_{k};
\end{equation}
and
\begin{equation}\label{}
    (\mu ^{L}_{ k})^{2}(x_{k})=\varphi ^{-a_{k}} (x_{k}) \ \ \text{for some nonnegative integer} \ \ a_{k}.
\end{equation}

\end{lem}
\begin{proof}  Since $\mu^{R}_{k}(\Gamma_{k})=-\Gamma_{k}$. Then,  one has
\begin{eqnarray*}
% \nonumber to remove numbering (before each equation)
  (\mu^{R}_{k})^{2}(x_{k}) &=& \mu^{R} _{k}(\Large(\prod_{i, \cdot_{i}\rightarrow \cdot_{k}} f_{ki}^{v_{ik}}+\prod_{i, \cdot_{k}\rightarrow \cdot_{i}} f_{ki}^{v_{ki}}\Large))x^{-1}_{k}) \\
  &=&  \Large(\prod_{i, \cdot_{i}\rightarrow \cdot_{k}} f_{ki}^{v_{ik}}+\prod_{i, \cdot_{k}\rightarrow \cdot_{i}} f_{ki}^{v_{ki}}\Large)x_{k}\Large(\prod_{i, \cdot_{i}\rightarrow \cdot_{k}} f_{ki}^{v_{ik}}+\prod_{i, \cdot_{k}\rightarrow \cdot_{i}} f_{ki}^{v_{ki}}\Large)^{-1} \\
      &=& \varphi ^{a_{k}} (x_{k}).
\end{eqnarray*}
The last equation is by the commutativity of $x_{1},\ldots, x_{k},\ldots x_{n}$ and applying (3.8), noticing that $\Gamma_{k}$ is a balanced star quiver with frozen rank $a_{k}$. This finishes the proof of equations (3.11). The proof of  (3.12) is quite similar  except for using the commutation relations (3.9) instead of (3.8).
\end{proof}

\begin{cor}

Let $p_{n}=(F, X, \Gamma, \varphi)$  be a  $\varphi$-commutative preseed with $\varphi$  is a finite order ring homomorphism. If $\Gamma_{k}$ is a balanced star quiver,  then  there is a non negative  integer $r$ such that

\begin{equation}\label{}
   (\mu^{R}_{k})^{2r}(p_{n} )=(\mu^{L}_{k})^{2r}(p_{n})=p_{n}.
\end{equation}

\end{cor}
\begin{proof}
 Assume that  $\varphi ^{r}=id_{\mathcal{D}_{n}}$ for some non negative integer $r$. Then using (3.11) ($r$-times) we get
\begin{equation}\label{}
    (\mu^{R}_{k})^{2r}(x_{k})=\varphi^{ra_{k}}(x_{k})=x_{k}.
\end{equation}
And, we already have  $(\mu^{R}_{k})^{2r}(x_{j})=x_{j}$ for $j\neq k$, and $(\mu^{R}_{k})^{2r}(\Gamma)=\Gamma$, which finishes the proof.
\end{proof}

\begin{que} For which  preseed  $p_{n}=(F, X, \Gamma)$, the set of cluster variables  $\mathcal{X}(p_{n})$ is finite?
\end{que}

In the following we  provide  equivalent conditions on $\varphi$-commutative  preseed to produce  a finite type cluster algebra.

\begin{thm} Let $p_{n}=(F, X, \Gamma, \varphi)$ be a balanced, $\varphi$-commutative  preseed. Then, the set of all cluster variables  $\mathcal{X}(p_{n})$ is a finite set if and only if $\varphi$ is of finite order.
\end{thm}

\begin{proof} Let $\mu_{i_{i}}, \ldots, \mu_{i_{t}}$ be  a sequence of mutations containing $j$ copies of $\mu_{k}$. Then,  by the definition of preseeds  mutation, we have
\begin{equation}\label{}
   \nonumber \mu_{i_{i}}\cdots \mu_{i_{t}}(x_{k})=\mu^{j}_{k}(x_{k}).
\end{equation}
Hence, for any preseed  $(F, X, \Gamma)$, we have

 \begin{equation}\label{}
   \mathcal{X}(p_{n})=\bigcup^{n}_{k=1}\mathcal{X}(p_{1}(x_{k})), \  \text{where} \ p_{1}(x_{k})=(\{F_{k}\}, \{x_{k}\}, \{\Gamma_{k}\}).
 \end{equation}

So, from (3.11) one has, for every $k \in [1, n]$ the set of all cluster variables of $p_{1}(x_{k})$ contains the set of the cluster variables of the form $\{\varphi^{la_{k}}(x_{k}); \ l\in \mathbb{N}\}$, where $a_{k}$ is the frozen rank of $p_{1}(x_{k})$. The set $\{\varphi^{la_{k}}(x_{k}); \  l\in \mathbb{N}\}$ is an infinite set if $\varphi$ is not of finite order. Which implies that if the set of cluster variables of $p_{n}$ is  finite  then $\varphi$ must be of finite order. Now assume that, $\varphi$ is of finite order. Then from (3.13), the seed $p_{1}(x_{k})$  will be reproduced after applying $\mu_{k}$,  $2r$-times which means that the set of  cluster variables of $p_{1}(x_{k})$  is finite for every $k\in [1, n]$ and then so is the set of cluster variables of $p_{n}$.
  \end{proof}

[16] Let $f$ be a $R_{n}$-linear automorphism over $\mathcal{D}_{n}$. Then $f$ is said to be a \emph{cluster variable preserver} of the preseed  $p_{n}$ if it keeps the cluster set $\mathcal{X}$ of $p_{n}$ invariant.

\begin{que} Given a preseed   $p_{n}=(F, X, \Gamma)$  describe the set of all cluster preservers of $p_{n}$.
\end{que}
\begin{thm} Let $p_{n}=(F, X, \Gamma, \varphi)$ be a balanced $\varphi$-commutative preseed  with frozen rank $(a_{1},\ldots, a_{n})$. Let $\phi_{l}$ be  the $R_{n}$-linear automorphisms of $\mathcal{D}_{n}$ induced by
\begin{equation}\label{}
    \phi_{l}(t)=t, \ \ \forall t\in R_{n} \ \ \text{and} \ \  \phi_{l}(x_{k})=\varphi^{la_{k}}(x_{k}), \ \ \forall k  \in [1,n].
\end{equation}
Then, for every $l\in \mathbb{Z}$, $\phi_{l}$ is a  cluster variables preserver for $p_{n}$.
\end{thm}

\begin{proof} Notice that, by definition of $\phi _{l}$,  it  depends  on  the frozen rank of $p_{n}$, which is invariant under mutation,  thanks to Proposition 3.8.

 First, for nonnegative integers.
 Let $l=1$. Equations (3.11) assure that, the action of the  automorphism $\phi _{1}$ on the cluster variables of $p_{n}$ is  identified with the action of the sequence of the mutations $\prod^{n} _{i=1}(\mu^{R}_{i})^{2}$.

     Let $x$ be an element of $\mathcal{X}(p_{n})$, without loss of generality, we  assume that   $x$ is a cluster variable of some seed  $q_{n}$, that can be obtain from $p_{n}$ by applying some sequence of only right mutations say $\mu^{R} _{i_{1}}\ldots \mu^{R}_{i_{d}}$. Then, $\phi_{1}(x)$ must be  a cluster variable in the seed $\prod^{n} _{i=1}(\mu^{R}_{i})^{2}(q_{n})=\prod^{n} _{i=1}(\mu_{i}^{R})^{2}\mu^{R} _{i_{1}}\ldots \mu^{R}_{i_{d}}(p_{n})$.\\ For $l\geq 2$,  again using (3.11), the action of   $\phi_{l}$ is identified with the action of the sequence of mutations  $(\prod^{n} _{i=1}(\mu^{R}_{i})^{2})^{l}(p_{n})$. Proving that $\phi_{l}$ permutes the elements of  $\mathcal{X}(p_{n})$ is quite similar to the case of $l=1$  with the obvious changes.\\ The case, when $l$ is a  negative integer, is similar, with using equations (3.12) instead of equations (3.11).
\end{proof}

\section{Weyl cluster algebras}

\subsection{Definition of generalized Weyl algebras}

\begin{defn}[Generalized Weyl algebra (1, 14, 15)] Let $\{\xi_{1},\ldots,\xi_{n}\}$ be a fixed set of elements of a commutative ring $R$ and $\theta =\{\theta _{1},\ldots, \theta _{n}\} $ be
a set of ring  automorphisms  such that $\theta_{i} (\xi_{j})=\xi_{j}$ for all $i\neq j$. \emph{The generalized Weyl  algebra} of degree $n$, denoted by $R(\theta, \xi,n)$, is defined to be the  ring extension of  $R$ generated by the  $2n$ indeterminates  $x_{1}, \ldots, x_{n}, y_{1}, \ldots,y _{n}$ modulo the commutation relations:
\begin{equation}\label{}
    x_{i}r=\theta_{i}(r)x_{i} \ \ \text {and} \ \ y_{i}r=\theta_{i} ^{-1}(r)y_{i}, \ \ \text {for any } i\in [1,n] \ \text {and for any}\ r\in R,
   \end{equation}
\begin{equation}\label{}
  x_{i}y_{i}=\xi_{i}, \ y_{i}x_{i}=\theta ^{-1}(\xi_{i}), \  x_{i}y_{j}=y_{j}x_{i}, \ \ x_{i}x_{j}=x_{j}x_{i} \ \text{and} \ y_{i}y_{j}=y_{j}y_{i}\ \ \forall i\neq j \in [1,n].
\end{equation}
We warn the reader that $x_{i}y_{i}\neq y_{i}x_{i}$ in general. The variables  $x_{1}, \ldots, x_{n}, y_{1}, \ldots,y _{n}$ are called \emph{Weyl variables}.
\end{defn}

\begin{exam}[4, 14, 15] Let $A_{n}$ be the $n^{th}$ Weyl algebra  generated by the $2n$ variables $x_{1}, \ldots, x_{n},y_{1}\ldots, y_{n}$ over  $K$ with the  relations
\begin{equation}\label{}
   x_{i}y_{i}-y_{i}x_{i}=1, \  \text{and} \ \ x_{i}x_{j}=x_{j}x_{i}, \ \ y_{i}y_{j}=y_{j}y_{i} \ \ \text{for} \ \ i \neq j, \ \forall i, j \in [1, n].
\end{equation}
 Let $\xi _{i}=y _{i}x _{i}+1$,  $R$ be the ring of polynomials $K[\xi _{1},\ldots, \xi_{n}]$ and $\theta_{i}: R\rightarrow R$, induced by $\xi _{i}\mapsto \xi_{i}+1, \xi_{j}\mapsto \xi_{j}, j\neq i,\ \text{for all} \ i,j\in[1, n]$. It is known that $A_{n}$ is isomorphic to the generalized Weyl algebra $R(\theta, \xi,n)$.
\end{exam}

\begin{exam}[14, 15] The coordinate algebra $A(SL_{q}(2,k))$ of algebraic quantum group $SL_{q}(2, k)$ is the $K$-algebra generated by  $x, y, u$, and $v$ subject to the following relations
\begin{equation}\label{}
   qux=xu, \ \ qvx=xv, \ \ qyu=uy, \ \ qyv=vy, \ \ uv=vu, \ \ q\in K^*
\end{equation}
\begin{equation}\label{}
   xy=quv+1, \ \  \text{and} \ \ yx=q^{-1}uv+1.
\end{equation}
$A(SL_{q}(2,k))$ is isomorphic to the generalized Weyl algebra $R(\xi, \theta, 1)$, where  $R$ is the algebra of polynomials $K[u,v]$; $\xi=1+q^{-1}uv$ and $\theta$ is an automorphism of $R$, defined by $\theta (f(u,v))=f(qu,qv)$ for any polynomial $f(u,v)$.

\end{exam}

\begin{defn} [Weyl preseeds and $q$-commutative preseeds] Let $p_{n}=(F, X, \Gamma)$ be a preseed  of rank $n$ in $\mathcal{D}_{n}$.  A quadruple $(F, X, \Gamma,  \theta)$  is said to be a \emph{Weyl preeseed } if there is a set $\theta= \{\theta_{1},\ldots, \theta_{n}\}$ of ring automorphisms  of $R_{n}$,  such that  for every $i \in [1, n]$,  $\theta_{i}$ fixes all the exchange cluster variables and  satisfies

 \begin{equation}\label{}
    x^{\pm 1}_{k}f=\theta_{i}^{\pm 1} (f)x^{\pm 1}_{i}, \ \ \forall f \in F_{i},  \forall i \in[1, n].
\end{equation}
If there is a fixed scalar $q\in K^{\ast}$ such that  $\theta _{i}$ satisfies

\begin{equation}\label{}
    \theta _{i}(f_{i})=qf_{i}, \ \  \text{for every}  \ \ \ i\in[1, n].
\end{equation}
In such special case,  $p_{n}=(F, X, \Gamma, q)$  is called  \emph{$q$-commutative preseed}.
\end{defn}

Let $p_{n}=(F, X, \Gamma)$ be a preseed and let \begin{equation}\label{}
    \xi_{k}=\prod_{ \cdot_{i}\rightarrow \cdot_{k}} f_{ki}^{v_{ik}}+\prod_{ \cdot_{k}\rightarrow \cdot_{i}} f_{ki}^{v_{ki}}, \ \ \  k\in [1, n].
\end{equation}
Then Relations 4.6 can be extended to $\xi_{k}^{\pm 1}$  as follows
\begin{equation}\label{}
  \theta_{k}^{\pm 1}(\xi_{k}^{\mp 1})x_{k}^{\pm 1}=x_{k}^{\pm 1}\xi_{k}^{\mp 1}, \ \ \ \  k\in [1, n].
\end{equation}

\begin{exam} Let $R(\theta, \xi, n)$ be a generalized  Weyl algebra.  Consider the quintuple   $p_{n}=(F, Y, \Gamma, \varphi, \theta )$, where $F=\{F_{i}\}^{n}_{i=1}$, $F_{i}=\{f_{i};  f_{i}=y_{i}x_{i}\}$, $Y=\{y_{1}, \ldots, y_{n}\}$, $\Gamma=\{\Gamma_{i}\}^{n}_{i=1}$ such that for $i\in [1, n], \Gamma_{i}$ is the  quiver

\begin{equation}\label{}
\nonumber     \xymatrix{
 \cdot_{f_{i}} &  \cdot_{y_{i}} \ar[l] },
\end{equation}
and $\varphi$ is given by

\begin{equation}\label{}
   \ \ \varphi(y_{i})=\xi_{i}y_{i}\xi^{-1}_{i}, \ \text{where} \       \xi_{i}=1+f_{i}, \  i \in [1, n].
\end{equation}
A short calculation shows that $\varphi$ satisfies Equations (3.7), hence  $p_{n}=(F, X, \Gamma, \varphi, \theta )$ is a $\varphi$-commutative   preseed. Also, from the properties of the $R$-automorphisms $\theta=(\theta_{1}, \ldots, \theta_{n})$  given in Equations (4.1) and (4.2)  one can see that $\theta_{i}$ satisfies Equation (4.6)   for each $i\in [1, n]$ which makes $p_{n}$ a Weyl preseed, then $p_{n}$ is $\varphi$-commutative Weyl  preseed.

In this case the iterated division rings   $\mathcal{D}_{i}, \ i=1, \ldots, n$, attached with $p_{n}$, are subrings of the division ring of rational functions in $y_{1}, \ldots, y_{n}$ over the ring $R$. In particular, in the case of  the $n^{th}$ Weyl algebra  $A_{n}$, the ring $R$ is ring of polynomials $K[\xi_{1}, \ldots, \xi_{n}]$. One can see that this ambient division ring of rational functions is an Ore domain. For information about Ore domains we refer to [13, 3].
\end{exam}

\begin{exam} Recall  the coordinate algebra $A(SL_{q}(2,k))$ of the algebraic quantum group $SL_{q}(2, k)$.  Consider the preseed  $p_{1}=(\{F_{1}\}, \{x\}, \{\Gamma_{1}\}, \{\theta_{1}\})$, where $F_{1}=\{qu, v\}$, $\theta_{1}: R\rightarrow R$ given by
 $\theta_{1} (f(u,v))=f(qu,qv)$  and  $\Gamma_{1} $ is given by
\begin{equation}\label{}
     \xymatrix{
 \cdot_{qu}    & \cdot_{v}\\
  \cdot_{x} \ar[u] \ar[ur]}.
\end{equation}
One can see that $p_{1}$ is a $q$-commutative preseed. Let $\zeta=quv+1$.  The cluster set of $p_{1}$ is given by

$\mathcal{X }(p_{1})=\{x,\zeta^{j}x\zeta^{-j}, \zeta ^{j+1}x^{-1}\zeta ^{-j-1}, j\in
\mathbb{N} \}\bigcup\{y,\zeta^{j}y\zeta^{-j}, \zeta ^{j+1}y^{-1}\zeta ^{-j-1},
j\in \mathbb{N}\}$.
\end{exam}

\begin{rem} If $p_{n}=(F, X, \Gamma,  \theta)$ is a  Weyl  preseed, then  the two quadruples $(F, \mu^{R}_{i}(X), -\Gamma, \theta^{-1})$ and $(F, \mu^{L}_{i}(X), -\Gamma, \theta^{-1})$ where $\theta=\{\theta_{1} ^{-1},\ldots, \theta_{n} ^{-1}\}$ are again   Weyl preseeds, for every $i\in [1, n]$.
\end{rem}

\begin{lem} Let $p_{n}$ be a    $q$-commutative preseed with  $q$ being an $m^{th}$ root of unity, for some natural number $m$ and let  $\Gamma_{k}$ be  balanced star quiver. Then, we have
\begin{equation}\label{}
   (\mu^{R}_{k})^{2m}(p_{n} )=(\mu^{L}_{k})^{2m}(p_{n})=p_{n}.
\end{equation}
\end{lem}
\begin{proof}

\begin{eqnarray*}
% \nonumber to remove numbering (before each equation)
  (\mu^{R}_{k})^{2}(x_{k}) &=& \mu^{R} _{k}(\Large(\prod_{i, \cdot_{i}\rightarrow \cdot_{k}} y_{i}^{v_{ik}}+\prod_{i, \cdot_{k}\rightarrow \cdot_{i}} y_{i}^{v_{ki}}\Large))x^{-1}_{k}) \\
  &=&  \Large(\prod_{i, \cdot_{i}\rightarrow \cdot_{k}} y_{i}^{v_{ik}}+\prod_{i, \cdot_{k}\rightarrow \cdot_{i}} y_{i}^{v_{ki}}\Large)x_{k}\Large(\prod_{i, \cdot_{i}\rightarrow \cdot_{k}} y_{i}^{v_{ik}}+\prod_{i, \cdot_{k}\rightarrow \cdot_{i}} y_{i}^{v_{ki}}\Large)^{-1} \\
   &=& \Large(\prod_{i, \cdot_{i}\rightarrow \cdot_{k}}y_{i}^{v_{ik}}+\prod_{i, \cdot_{k}\rightarrow \cdot_{i}} y_{i}^{v_{ki}}\Large)x_{k}\Large(\prod_{i, \cdot_{i}\rightarrow \cdot_{k}} y_{i}^{v_{ik}}+\prod_{i, \cdot_{k}\rightarrow \cdot_{i}} y_{i}^{v_{ki}}\Large)^{-1}\\
   &=& \Large(\prod_{i, \cdot_{i}\rightarrow \cdot_{k}}y_{i}^{v_{ik}}+\prod_{i, \cdot_{k}\rightarrow \cdot_{i}} y_{i}^{v_{ki}}\Large)\theta_{k}(\Large(\prod_{i, \cdot_{i}\rightarrow \cdot_{k}} y_{i}^{v_{ik}}+\prod_{i, \cdot_{k}\rightarrow \cdot_{i}} y_{i}^{v_{ki}}\Large)^{-1})x_{k}\\
    &=& \Large(\prod_{i, \cdot_{i}\rightarrow \cdot_{k}}y_{i}^{v_{ik}}+\prod_{i, \cdot_{k}\rightarrow \cdot_{i}} y_{i}^{v_{ki}}\Large)(\theta_{k}(\Large\prod_{i, \cdot_{i}\rightarrow \cdot_{k}} y_{i}^{v_{ik}}+\prod_{(i, \cdot_{k}\rightarrow \cdot_{i}} y_{i}^{v_{ki}}\Large))^{-1}x_{k}\\
   &=& q^{-a_{k}}x_{k}, \ \text{where $a_{k}$ is the frozen component of $\Gamma_{k}$}.
\end{eqnarray*}
Then
\begin{equation}\label{}
   \nonumber (\mu^{R}_{k})^{2m}(x_{k})=q^{-ma_{k}}x_{k}=x_{k}.
\end{equation}
And since we already have $(\mu^{R}_{k})^{2m}(\Gamma)=\Gamma$, which completes the proof.
\end{proof}

\begin{cor} Let $p_{n}=(F, X, \Gamma, q)$ be a $q$-commutative balanced preseed. Then the set of of all cluster variables $\mathcal{X}(p_{n})$ is finite if and only if $q$ is an $m^{th}$-root of unity, for some natural number $m$.
\end{cor}
\begin{proof}
  From Lemma 4.8 and Equation (3.15) we have

  \begin{equation}\label{}
\nonumber  \{q^{(-a_{k})l}x_{k}; \  l\in \mathbb{N}\}  \subset\mathcal{X}(p_{1}(x_{k}))\subset \mathcal{X}(p_{n}),  \forall k \in [1, n].
  \end{equation}
 If for every natural number $m$, $q^{m}\neq 1$,  then for each $k\in[1, n]$,  the set  $\{q^{(-a_{k})l}x_{k}; \  l\in \mathbb{N}\}$ is an infinite set. So, if $\mathcal{X}(p_{n})$ is a finite set then $q$ must be an $m^{th}$-root of unity, for some natural number $m$. Now, assume that $q^{m}= 1$, for some natural number $m$, then again using Lemma 4.8,  for each $k$ in [1, n], the seed $p_{1}(x_{k})$ will be reproduced after applying $\mu^{R}_{k}$, $2m$-times. Then $\mathcal{X}(p_{1}(x_{k}))$ is a finite set for each $k$ in $[1, n]$ and hence again from (3.15), the set $\mathcal{X}(p_{n})$ must be a finite set.

\end{proof}

\begin{defn}[Weyl cluster algebras]  Let  $p_{n}=(F, X, \Gamma, \theta)$  be a Weyl preseed. The Weyl cluster algebra $\mathcal{H}(p_{n})$  is defined to be the  $R_{n}$-subalgebra of $\mathcal{D}_{n}$  generated by the cluster set  $\mathcal{X} (p_{n})$.
\end{defn}

The following remark and theorem shed some light on the structure of the Weyl cluster algebra $\mathcal{H}(p_{n})$. Remark 4.11 and first part of the Theorem 4.12 can be phrased as, the Weyl cluster algebra $\mathcal{H}(p_{n})$ is generated by $R_{n}$ and many (could be infinite) isomorphic copies of generalized Weyl algebras, each vertex in the exchange graph of $p_{n}$ gives rise to two copies of them. The second part of the theorem is the Laurent phenomenon, Theorem 2.8, in the  Weyl preseeds case.  The third part of the same theorem is simply saying that  $\mathcal{H}(p_{n})$ is isomorphic to the tensor product of  the $n$  Weyl cluster algebras of rank one associated to the $n$ iterated rank one Weyl preseeds  associated to $p_{n}$.

\begin{rem and def} Let $p_{n}=(F, X ,\Gamma, \theta)$ be a Weyl preseed and $R=K[\xi_{1}, \ldots, \xi_{n}]$ be the ring of polynomials in $\xi_{1}, \ldots, \xi_{n}$ where $\xi_{i}, i=1,\ldots, n$ are as defined in (4.8). Then $p_{n}$ gives rise to two copies of generalized Weyl  algebras of rank $n$, as follows
\begin{enumerate}
  \item [(a)] $H^{R}(p_{n})$ is the ring extension of $R$ generated by $\mu^{R}_{1}(x_{1}), \ldots, \mu^{R}_{n}(x_{n}), x_{1}, \ldots, x_{n}$.
  \item [(b)] $H^{L}(p_{n})$ is the ring extension of $R$ generated by $x_{1}, \ldots, x_{n}, \mu^{L}_{1}(x_{1}), \ldots, \mu^{L}_{n}(x_{n})$.
  \item [(b)] In particular, if $p_{n}=(F, Y ,\Gamma, \varphi, \theta)$ is the preseed given in Example 4.5, then each of $H^{R}(p_{n})$ and $H^{L}(p_{n})$ are isomorphic to $R(\theta, \xi, n)$ as generalized Weyl algebras. In the case of $H^{R}(p_{n})$  (respectively $H^{L}(p_{n})$) the isomorphism is defined by  sending the cluster variable $\mu^{R}_{i}(x_{i})$ to Weyl variable  $y_{i}$ and the cluster variable $x_{i}$ to the Weyl variable $y_{i}$ of $R(\theta, \xi, n)$  (respectively by  sending the cluster variable $x_{i}$ to the Weyl variable $x_{i}$ and $\mu^{R}_{i}(x_{i})$ to the Weyl variable $y_{i}$) for $i=1,\ldots, n$. Details for the case $n=1$ are given in Example 4.14.
\end{enumerate}

\end{rem and def}

\begin{thm} Let $p_{n}=(F, X ,\Gamma, \theta)$ be a Weyl preseed in $\mathcal{D}_{n}$. Then the following are true
\begin{enumerate}

\item Right and left mutations on $p_{n}$ induce isomorphisms between the generalized Weyl algebras $H^{R}(p_{n})$ and $H^{R}(\mu^{R}_{k}(p_{n}))$ (respectively $H^{L}(p_{n})$ and $H^{L}(\mu^{L}_{k}(p_{n}))$).
\item The  Weyl cluster algebra $\mathcal{H}(p_{n})$ is a subring of the (non-commutative)  ring of Laurent polynomials in the initial exchange cluster variables with coefficients from ring of polynomials $R_{n}[\theta^{\pm 1}_{1}(\xi^{-1}_{1}), \ldots,  \theta^{\pm 1}_{n}(\xi^{-1}_{n})]$.

\item Let $p_{1}(x_{k})$ be the rank one preseed $(F_{k}, \{x_{k}\}, \{\Gamma_{k}\}, \theta_{k} \})$. Then

\begin{equation}\label{}
    \mathcal{H}(p_{n})\cong \mathcal{H}(p_{1}(x_{1}))\otimes \cdots \otimes \mathcal{H}(p_{1}(x_{n})).
\end{equation}

\end{enumerate}
\end{thm}
\begin{proof}

To prove part (1), consider the  $R_{n}$-linear automorphism of  $\mathcal{D}_{n}$, denoted by  $T^{R}_{p_{n},k}: \mathcal{D}_{n}\rightarrow \mathcal{D}_{n}$ induced by $x_{k}\mapsto \mu^{R}_{k}(x_{k})$, $k\in[1,n]$. The restriction of this automorphism on   $H^{R}(p_{n})$  induces the  algebras isomorphism   $\widehat{T}^{R}_{p_{n}, k}: H^{R}(p_{n}) \rightarrow H^{R}(\mu^{R}_{k}(p_{n}))$ given by $r\mapsto r, \forall r \in \mathcal{R}$, and  $x_{k}\mapsto \mu^{R}_{k}(x_{k})=\xi_{k} x^{-1}_{k}, \forall k \in [1,n]$. Which implies $\mu^{R}_{k}(x_{k})\mapsto \xi_{k} x_{k}
    \xi^{-1}_{k}=\mu^{R}_{k}(\mu^{R}_{k}(x_{k}))$. Finally, it is easy to see  that the generalized Weyl commutation relations (4.1)  are invariant under $\widehat{T}^{R}_{p_{n}, k}$. (The argument
    for $H^{L}(p_{n})$ is  quite similar).\\
For part (2), let $y\in \mathcal{X}(p_{n})$. Without loss of generality, using (3.15) we can assume that $y$ is an element of $\mathcal{X}(p_{1}(x_{k}))$ for some $k\in [1, n]$.  Hence,  $y$ can be obtained from $x_{k}$ by applying some sequence of mutations on $p_{1}(x_{k})$. Let $l$ be the length  of a shortest such sequence of mutations. By (3.6) we have that every non-trivial sequence of mutations can be reduced to either only right mutations or only left mutations. Then, by mathematical induction on $l$, one can show

\begin{equation}\label{}
   y=\left\{
                  \begin{array}{ll}
                    \xi_{k}^{\frac{l+1}{2}}x_{k}^{-1}\xi_{k}^{-(\frac{l+1}{2}-1)} \ \text{or} \ \ \xi_{k}^{-(\frac{l+1}{2}-1)}x_{k}^{-1}\xi_{k}^{\frac{l+1}{2}}, \  \text{if} \ l  \ \text{is an odd number;}& \\
                     \xi_{k}^{\frac{l}{2}}x_{k}\xi_{k}^{-\frac{l}{2}} \ \text{or} \  \  \xi_{k}^{-\frac{l}{2}}x_{k}\xi_{k}^{\frac{l}{2}}, \ \   \text{if} \ l \ \text{is an even number.}&
                  \end{array}
                \right.
\end{equation}

Now, let   $\textsl{m}$ be a monomial in the elements of  $\mathcal{X}(p_{n})$. Then again using (3.15), and the identities  (4.14), (4.9) and the commutations relations (4.2), one can write $m$ as  $rm'$ where $r$ is a monomial in the elements from the set $F^{n}\bigcup \{\theta^{\pm 1}_{1}(\xi^{-1}_{1}), \ldots,  \theta^{\pm 1}_{n}(\xi^{-1}_{n})\}$ and $m'$ is a monomial of elements from $\{x^{\pm 1}_{k}, \ k\in [1, n]\}$. Finally, the elements of $\mathcal{H}(p_{n})$ are finite sum of finite product of monomials from the elements of $\mathcal{X}(p_{n})$ which finishes the proof of Part (2).\\
For Part (3), by the definition of Weyl cluster algebra and the proof of Part (2) above  one can see that the Weyl cluster algebra $\mathcal{H}(p_{1}(x_{k}))$ is generated as a $K$-vector space by the monomials

\begin{equation}\label{}
   \mathfrak{m}_{k}=\{f^{\alpha_{k1}}_{k1}\cdots f_{km_{k}}^{\alpha_{km_{k}}}(\theta^{\pm 1}_{k}(\xi^{-1}_{k}))^{\alpha'_{k}}x^{\beta}_{k};   \alpha_{kj}, \alpha'_{k},   \beta\in \mathbb{Z}, \forall j\in [1, m_{k}]\}.
\end{equation}
 Then, the  Weyl cluster algebra $\mathcal{H}(p_{n})$ is generated as a vector space by $\mathfrak{m}(p_{n})=\{m_{1}\cdots m_{n}; \ m_{k}\in \mathfrak{m}_{k}, k\in [1, n]\}$, where $m_{i}m_{j}=m_{j}m_{i}$  for $m_{i}\in\mathfrak{m}_{i}$ and $m_{j}\in\mathfrak{m}_{j}$ for every  $i\neq j \in [1, n]$. Now we will show that $\mathfrak{m}(p_{n})$ consists of linearly independent elements. Let $E=R_{n}[\theta^{\pm 1}_{1}(\xi^{-1}_{1}), \ldots,  \theta^{\pm 1}_{n}(\xi^{-1}_{n})][t^{\pm 1}_{1}, \cdots, t^{\pm 1}_{n}]$, consider the linear endomorphisms of $E$  given by $X_{k}^{\pm 1}(f)=t^{\pm 1}f, f\in E$. The map $\sigma: \mathcal{H}(p_{n})\longrightarrow End(E)$ induced by $x_{k}^{\pm 1}\longmapsto X_{k}^{\pm 1}, k\in [1, n]$ defines an algebras homomorphism. One can see that the endomorphisms

\begin{equation}\label{}
    \nonumber f^{\alpha_{11}}_{11}\cdots f_{1m_{1}}^{\alpha_{1m_{1}}}\cdots f^{\alpha_{n1}}_{n1}\cdots f_{nm_{j}}^{\alpha_{nm_{n}}}(\theta^{\pm 1}_{1}(\xi^{-1}_{1}))^{\alpha'_{k1}} \ldots  (\theta^{\pm 1}_{n}(\xi^{-1}_{n}))^{\alpha'_{km_{k}}}X^{\beta_{1}}_{1}\cdots X^{\beta_{n}}_{n};\end{equation}
\begin{equation}\label{}
    \nonumber  \alpha_{ji}, \alpha'_{ji},  \beta_{j}\in \mathbb{Z}, i\in [1, m_{j}],  j\in [1, n]
\end{equation}
are linearly independent elements of $End(E)$ over $K$. Hence, $\mathfrak{m}(p_{n})$ consists of linearly independent elements which makes it  a basis for $\mathcal{H}(p_{n})$ as a $K$-vector space and $\sigma$ is an injective algebra homomorphism. Then the map that sends $m_{1}\cdots m_{n}$ onto $m_{1}\otimes\cdots \otimes m_{n}, m_{k}\in \mathfrak{m}_{k}, k\in [1, n]$ defines an isomorphism from  $\mathcal{H}(p_{n})$ to $\mathcal{H}(p_{1}(x_{1}))\otimes \cdots \otimes \mathcal{H}(p_{1}(x_{n}))$.

\end{proof}
From the Proof of Part (2) of Theorem 4.12, we have the following remark.
\begin{rem} The Weyl cluster algebra $\mathcal{H}(p_{n})$ is finitely generated algebra.
\end{rem}

\begin{exam}[Weyl cluster algebra associated to first Weyl algebra] Recall the $n^{th}$ Weyl algebra given in Example 4.2 and the associated preseed  given in Example 4.5. Let $A_{1}$ be the first Weyl algebra and consider the preseed
$p_{1}=(\{f\}, \{y\}, \{\xymatrix{ \cdot_{y} \ar[r] &\cdot_{f}}\})$. Here $R_{1}=K [\mathbb{P}]$, where $\mathbb {P}$ is the cyclic group generated by $f=yx$. Then We have the following exchange graph
\begin {itemize}
\item $\mathbb{G}(p_{1})$
\begin{equation}\label{}
    \xymatrix{ \ldots  \ar[r]_{L}&\ar[r]^{R} \stackrel{y_{-3}}{\cdot} \ar[l]_{R}&\ar[r]^{R} \stackrel{y_{-2}}{\cdot} \ar[l]^L &\ar[r]^{R} \stackrel{y_{-1}}{\cdot}\ar[l]^{L} &\ar[r]_{L} \stackrel{y_{0}=y}{\cdot}\ar[l]^{L} &\ar[l]_{R}\stackrel{y_{1}}{\cdot}\ar[r]^{R} &\ar[l]^L\stackrel{y_{2}}{\cdot}\ar[r]^R&\ar[l]^L\stackrel{y_{3}}{\cdot}\ar[r]^R & \ar[l]^{L} \ldots
    },
\end{equation}
\end {itemize}
(here $ \xymatrix{ & {\cdot}\ar[l]^{L}}$ is left mutation and $
\xymatrix{{\cdot}\ar[r]^R &}$ is right mutation). Which can be encoded by
the following equations
\begin{equation}\label{}
 y_{2k}y_{2k\pm 1}=y_{2k\pm 1}y_{2k}+1, \ \ \ \text{for} \ \  k\in \mathbb{Z}.
\end{equation}
The  Weyl cluster algebra  $\mathcal{H}(p_{1}(y))$ is the $R_{1}$-subalgebra of $\mathcal{D}_{1}$ generated by the set of cluster variables $\{y_{k}, k\in \mathbb{Z}\}$. Relations (4.17) can be interpreted as follows,  each arrow in $\mathbb{G}(p_{1})$ corresponds to a copy of  first Weyl algebra, denoted by $A^{k}_{1}=K\langle y_{k},y_{k+1}\rangle, k\in \mathbb{Z}$ and right (respectively left) mutations define isomorphisms between the adjacent copies, given by $T_{k}:A^{k}_{1}\rightarrow A^{k+1}_{1}, \  \ y_{k}\mapsto y_{k+1} $ for $k\in \mathbb{Z}$ (respectively to the inverses of $T_{k}, k\in \mathbb{Z}$).
\end{exam}

 The \emph{adjunction isomorphism}   $\Theta: R(\theta^{-1}, \theta^{-2}(\xi), 1)\rightarrow R(\theta, \xi, 1)$ given by $r\mapsto \theta^{-1}(r), x\mapsto y$ and $y\mapsto x$. In  [14], the adjunction isomorphism  played  an important role in describing the representations theory of generalized Weyl algebra $R(\theta, \xi)$.

\begin{rem} Consider the  preseed $p_{1}=(F, Y, \Gamma)$ associated to the generalized Weyl algebra $R(\xi, \theta, 1)$, given in Example 4.5.  The action of the adjunction isomorphism $\Theta$ on the exchange cluster variables of any  two adjacent seeds on the  exchange graph of $p_{1}$  coincides with the action of the right and left mutations.

\end{rem}

\section{Representations arising from Weyl cluster structure}
\subsection{Space of representations $V_{n}$}
 In the following, let $p_{n}=(F, Y, \Gamma, \theta)$ be the generalized Weyl preseed associated to the generalized Weyl algbera $R(\theta, \xi, n)$, as given in Example 4.5. A \textit{cluster monomial} in $\mathcal{H}(p_{n})$ is a product  of non negative powers of exchange cluster variables belonging to the same cluster. To visualize that, the monomial $m=z^{\beta _{1}}_{1}\cdots z^{\beta _{n}}_{n}, \beta _{i}\in \mathbb{Z}_{\geq 0}, i\in [1,n]$ is a cluster monomial if $\{z_{1},\ldots,z_{n}\}$ is the set of the exchange  cluster variables of some seed in the exchange graph of $p_{n}$.

\begin{defn} The \emph{space of representations} $V_{n}$ is defined to be the $K(f_{1},\ldots, f_{n})$-left span by the set of all cluster monomials.
\end{defn}

\begin{lem} The space of representations $V_{n}$ is independent of $p_{n}$ and  depends only on the exchange graph $\mathbb{G}(p_{n})$.
\end{lem}
\begin{proof} The statement of the lemma is equivalent to the fact that ``the set of all cluster monomials of every seed in $\mathbb{G}(p_{n})$ is the same" which is equivalent to ``any two  seeds  in $\mathbb{G}(p_{n})$ have the same exchange graph" which is an immediate result of the fact that the set of all seeds in $\mathbb{G}(p_{n})$ form an equivalent class under (left and right) mutations as equivalent relation which is due to Part (2) of  Proposition 3.4.
\end{proof}
\begin{prop} If $p_{n}$ is a  preseed, then the following are true
 \begin{enumerate}
 \item
For any set of $n$ (or less) different cluster variables, not including two variables  produced from the same initial cluster variable, there is at least one seed in $\mathbb{G}(p_{n})$  which contains all of them;
\item
For any two cluster variables $z_{1}$ and $z_{2}$, produced from the same initial cluster variable, there are two cases for their product
 \begin{itemize}
 \item if $z_{2}$ can be obtained from $z_{1}$ by applying  some sequence of mutations of an odd length, then $z_{1}z_{2} \in K(f_{1},\ldots, f_{n})$;
 \item if $z_{2}$ can be obtained from $z_{1}$ by applying some sequence of mutations of an even length, then $z_{1}z_{2}$ can be written as $g z^{2}_{1}$, for some $g\in K(f_{1},\ldots, f_{n})$.
 \end{itemize}
 \end{enumerate}
\end{prop}
\begin{proof} Every cluster variable can be traced back to one of the initial cluster variables. More precisely, for any $y\in \mathcal{X}(p_{n})$ there is $k\in [1, n]$ such that $y\in \mathcal{X}(p_{1}(x_{k}))$, thanks to (3.15). Hence, there is a sequence of mutations $\mu^{y}$ such that $y=\mu^{y}(x_{k})$. Now, let $\{y_{1},\ldots, y_{t}\}$ be a subset of $\mathcal{X}(p_{n})$ such that $t\in[1, n]$. Then, one can see that the cluster of the  seed $\mu^{y_{1}}\cdots\mu^{y_{t}}(p_{n})$ contains the set $\{y_{1},\ldots, y_{t}\}$.  Part (2) is immediate from (4.14).
\end{proof}
Let $Y=(y_{1},\ldots,y_{n})$ be the  cluster of the preseed  $p_{n}$. For $t\in\mathbb{Z}$, $y_{i, t}$ denotes the cluster variable obtained from the initial cluster variable $y_{i}$ by applying one of the following sequence of mutations   $(\mu^{R}_{i})^{t}$ if $t\geq 0$ or  $(\mu^{L}_{i})^{-t}$ if $t<0$.\\
Using Proposition 5.3 and the above notation, a typical element of $V_{n}$ can be written as  a sum of   elements  of the form
 \begin{equation}\label{}
    v=r(f_{1},\ldots, f _{n}) y^{\beta_{1}}_{1, m_{1}}\cdots y^{\beta_{n}}_{n, m_{n}},
 \end{equation}
 where  $r(f_{1},\ldots, f _{n})\in K(f_{1},\ldots, f_{n})$,  $(\beta_{1}, \ldots, \beta_{n})\in \mathbb{Z}_{\geq 0}^{n}$ and $(m_{1},\ldots, m_{n}) \in \mathbb{Z}^{n}$.

\begin{exam} Consider the Weyl preseed   $p_{n}=(F, Y, \Gamma, \theta)$, as given in Example 4.5. The $i^{th}$ branch of the exchange graph $\mathbb{G}(p_{n})$  is as follows

\begin{equation}\label{}
 \nonumber   \xymatrix{   \ldots \ar[r]^R \stackrel{(y_{1, m_{1}},\ldots, y_{i,m_{i}-1},\ldots , y_{n,m_{n}})}{\cdot}  &\ar[r]^R \stackrel{(y_{1,m_{1}},\ldots,y_{i,m_{i}},\ldots,y_{n,m_{n}})}{\cdot} \ar[l]^{L} & \stackrel{(y_{1,m_{1}}, \ldots,y_{i,m_{i}+1}, \ldots, y_{n,m_{n}})}{\cdot } \ar[l]^{L}\cdots}
\end{equation}
 For the sake of  simplicity, we labeled each vertex by the clusters only. The space of representations  $V_{n}$ is the left $K(\xi_{1},\ldots,\xi_{n})$-linear span  by the set
\begin{equation}\label{}
   \{y^{\beta _{1}}_{1,m_{1}} \cdots y^{\beta _{n}}_{n,m_{n}}| \ \  \text {for}  \ \ m=(m_{1}, \ldots, m_{n})\in \mathbb{Z}^{n}, \text {and}  \ \ \beta = (\beta _{1}, \ldots, \beta _{n})\in \mathbb{Z}^{n}_{\geq 0}\}.
\end{equation}
\end{exam}

\begin{defn}[Representations of $R(\theta, \xi, n)$  on $V_{n}$] An action of the generators $x_{1},\ldots, x_{n}$ and $y_{1},\ldots, y_{n}$  on the (a generic) element  $v$  (given in (5.1)), is given by
\begin{equation}\label{}
    y_{i}(v):= r(f_{1},\ldots, f_{i-1},\theta^{-1}_{i} (f_{i}),\ldots, f _{n})  y^{\beta_{1}}_{1,m_{1}}\cdots y^{\beta_{i-1}}_{i-1,m_{i-1}} y^{\beta_{i}}_{i,m_{i}-1} y^{\beta_{i+1}}_{i+1,m_{i+1}}\cdots y^{\beta_{n}}_{n,m_{n}};
\end{equation}
and
\begin{equation}\label{}
    x_{i}(v):=\theta_{i} (f_{i}) r(f_{1},\ldots, f_{i-1},\theta_{i} (f_{i}),\ldots, f _{n})  y^{\beta_{1}}_{1,m_{1}}\cdots y^{\beta_{i-1}}_{i-1,m_{i-1}} y^{\beta_{i}}_{i,m_{i}+1} y^{\beta_{i+1}}_{i+1,m_{i+1}}\cdots y^{\beta_{n}}_{n,m_{n}}.
\end{equation}

\end{defn}

\begin{lem}
 The actions given in definition 5.5 define  a fully faithful left module structure of  $R(\theta, \xi, n)$ on $V_{n}$.
\end{lem}

\begin{proof}
The module structure of $R(\theta, \xi, n)$ on $V_{n}$ is defined by extending  (5.3) and (5.4) to random elements of $R(\theta, \xi, n)$.
 It is obvious to see that the actions given in (5.3) and (5.4) are compatible with Relations (4.1). In the following we show that Relations (4.2) are satisfied on the generic element $v$, given in (5.1). We have

\begin{eqnarray}
  % \nonumber to remove numbering (before each equation)
 \nonumber x_{i}y_{i}(v)&=& x_{i}(r(f_{1},\ldots, \theta^{-1}_{i} (f_{i}),\ldots, f _{n})  y^{\beta_{1}}_{1,m_{1}}\cdots y^{\beta_{i-1}}_{i-1,m_{i-1}} y^{\beta_{i}}_{i,m_{i}-1} y^{\beta_{i+1}}_{i+1,m_{i}+1}\cdots y^{\beta_{n}}_{n, m_{n}}))\\
 \nonumber &=& \theta_{i} (f) r(f_{1},\ldots, \theta^{-1}_{i} (\theta (f_{i})), \ldots, f_{n} )  y^{\beta_{1}}_{1,m_{1}}\cdots y^{\beta_{i-1}}_{i-1,m_{i-1}} y^{\beta_{i}}_{i,m_{i}} y^{\beta_{i+1}}_{i+1,m_{i}+1}\cdots y^{\beta_{n}}_{n,m_{n}})\\
  \nonumber &=& \theta_{i}(f_{i})(v)\\
  \nonumber &=& \xi_{i} v.
  \end{eqnarray}
  And
  \begin{eqnarray}
    % \nonumber to remove numbering (before each equation)
    \nonumber y_{i}x_{i}(v) &=& y_{i}(\theta_{i} (f_{i}) r(f_{1},\ldots,\theta_{i} (f_{i}),\ldots, f _{n})  y^{\beta_{1}}_{1,m_{1}}\cdots y^{\beta_{i-1}}_{i-1,m_{i-1}} y^{\beta_{i}}_{i,m_{i}+1} y^{\beta_{i+1}}_{i+1,m_{i}+1}\cdots y^{\beta_{n}}_{n,m_{n}})) \\
    \nonumber &=&  \theta_{i} (\theta_{i} ^{-1}(f_{i})) r(f_{1},\ldots, \theta_{i} (\theta ^{-1}(f_{i})),\ldots, f _{n})  y^{\beta_{1}}_{1,m_{1}}\cdots y^{\beta_{i-1}}_{i-1,m_{i-1}} y^{\beta_{i}}_{i,m_{i}} y^{\beta_{i+1}}_{i+1,m_{i}+1}\cdots y^{\beta_{n}}_{n,m_{n}})\\
    \nonumber &=& f_{i}v.
  \end{eqnarray}
In a similar way, one can get the rest of the Relations (4.2). The property of fully faithful module is a straightforward from the definitions of the actions given in (5.3) and (5.4).

\end{proof}

\begin{prop}
 The  module structure given in Definition 5.5 can be extended to the Weyl cluster algebra associated to $p_{n}$.
\end{prop}
\begin{proof}  To upgrade the representations of $R(\theta, \xi, n)$ on $V_{n}$ to the Weyl cluster algebra  associated to $p_{n}$, we introduce the action of $y^{-1}_{i}$ on the element $v$, by $y^{-1}_{i}(v)=\theta^{-1}(\xi_{i})x_{i}(v)$. The action of a random element of the Weyl cluster algebra $\mathcal{H}(p_{n})$ will be induced from the action of both of $y_{i}$ and $y^{-1}_{i}$ for $i=1,\ldots, n$, thanks to Part (2) of Theorem 4.12.
\end{proof}

\subsection{Cluster strands and the strand submodules of $V_{n}$.} Before introducing the cluster strands we need to introduce the following notations. For $t\in \mathbb{Z}$, let

\begin{equation}\label{}
    \nonumber \theta ^{t}(-)= \begin{cases} \overbrace{\theta (\theta (\ldots \theta (-)))}^{t-times}, & \text{ if } t > 0 ,\\
         id_{R}, & \text{ if } t=0,\\
         \overbrace{\theta ^{-1} (\theta  ^{-1} (\ldots \theta  ^{-1} (-)))}^{|t|-times}, & \text{ if } t < 0.
\end{cases}
\end{equation}
Consider the following three sets of monomials in the elements
$\{\theta ^{t}(z); t \in \mathbb{Z} \}$

\begin{enumerate}
\item
 \begin{equation}\label{}
     \nonumber M^{+}(z ):=\{ 1, \theta ^{t}(z^{\pm 1})\theta ^{t+1}(z^{\pm 1})\cdots \theta ^{t+q}(z^{\pm 1})| q, t\in \mathbb{Z}_{\geq 0} \};
 \end{equation}
\item
   \begin{equation}\label{}
     \nonumber M^{-}(z ):=\{ 1, \theta ^{t}(z^{\pm 1})\theta ^{t-1}(z^{\pm 1})\cdots \theta ^{t-q}(z^{\pm 1})| q \in  \mathbb{Z}_{\geq 0}, t\in \mathbb{Z}_{<0} \};
 \end{equation}
\item\begin{equation}\label{}
    M(z):=\{ \emph{m}_{1}\emph{m}_{2}| \emph{m}_{1} \in M^{+}(z ) \ \  \text {and}  \ \  \emph{m}_{2} \in M^{-}(z )\}.
 \end{equation}
\end{enumerate}
 For every $h\in K(f_{1},\ldots, f_{n})$ and $t=(t_{1},\ldots, t_{n})\in \mathbb{Z}^{n}$  we associate a subset of $K(f_{1},\ldots, f_{n})$ as follows
 \begin{equation}\label{}
 c(h, t):=\{\alpha_{1}\cdots\alpha_{n}  h(\theta^{t _{1}} _{1}(f_{1}),\ldots,\theta^{t _{n}} _{n}(f_{n}))| \  \alpha_{i}  \in  M(f_{i}), \forall i\in[1,n] \}.
\end{equation}

\begin{defn}[\textit{Cluster strands}] Fix a natural number  $l$ and  a one to one map  $\sigma: [1, l]\rightarrow \mathbb{Z}^{n}_{\geq 0}\times \mathbb{Z}^{n}$. Let  $\beta=(\beta _{1},\ldots, \beta _{l}) \in \overbrace{\mathbb{Z}^{n}_{\geq 0}\times \cdots \times \mathbb{Z}^{n}_{\geq 0}}^{l-\text{times}}$ and $m =(m _{1},\ldots, m _{l}) \in \overbrace{\mathbb{Z}^{n}\times \cdots \times\mathbb{Z}^{n}}^{l-\text{times}}$ such that $\sigma (j)=(\sigma _{1}(j),\sigma _{2}(j))=(\beta_{j}, m _{j})$   where $\beta_{j}=(\beta_{j1},\ldots, \beta_{jn})$ and $m _{j}=(m_{j1},\ldots, m_{jn}), j\in [1,l]$. Let   $r=(r_{1},\ldots ,r_{l})$ such that $r_{j}\in K(f_{1},\ldots, f_{n})$ for  $j\in[1,l]$. Consider the following subset of  $V_{n}$

\begin{equation}\label{}
  S_{l}(\sigma, r):=\Big\{\sum^{l}_{j=1}g_{j}y^{\beta_{j1}}_{1,m_{j1}+t_{j1}}\cdots y^{\beta_{jn}}_{n,m_{jn}+t_{jn}}| g_{j}\in c(r_{j}, t_{j}), t_{j}=(t_{j1},\ldots, t_{jn})\in \mathbb{Z}^{n},  j \in [1, l]\Big\}.
\end{equation}
With the above data, $S_{l}(\sigma, r)$ is called a \textit{cluster strand} \emph{of length $l$, with respect to $r$  and $\sigma$}. Furthermore, $S_{l}(\sigma, r)$ is called a \emph{full cluster strand} if $\sigma _{1}(j)\in \mathbb{Z}^{n}_{>0}$ for every $j\in [1, l]$.
\end{defn}
\begin{exam} [A cluster strands of  length $2$ in $V_{3}$] Let $l=2$, $\sigma_{1} (1)=(0,3,0), \sigma_{1} (2)=(1, 0,2)$, $\sigma _{2}(1)=(1, 1, 0), \sigma _{2}(2)=(0,1, 1)$, and $r=(f^{2}_{1}+f_{2}, f_{1}f_{3})$. For $t_{j}=(t_{j1}, t_{j2}, t_{j3})\in \mathbb{Z}^{3}, j\in [1, 2]$, we have
\begin{equation}\label{}
  \nonumber c(f^{2}_{1}+f_{2},  t_{1})=\{ \alpha_{1}\alpha_{2} \alpha_{3} ((\theta^{t _{11}} _{1}(f_{1}))^{2}+\theta^{t _{12}} _{2}(f_{2}))| \  \alpha_{i}  \in  M(f_{i}),  i\in[1,3] \},
\end{equation}
and
\begin{equation}\label{}
  \nonumber c( f_{1}f_{3}, t_{2})=\{ \alpha_{1}\alpha_{2}  \alpha_{3}\theta^{t _{21}} _{1}(f_{1})\theta_{3}^{t_{23}} (f_{3})| \  \alpha_{i}  \in  M(f_{i}),  i\in[1,3] \}.
\end{equation}
With the above data we have
\begin{equation}\label{}
   \nonumber S_{3}(\sigma, r)=\big\{g_{1}y^{3}_{2,1+t_{12}}+g_{2} y_{1,0+t_{21}}y^{2}_{3,1+{t_{23}}}|g_{1}\in c(f^{2}_{1}+f_{2}, t_{1}),g_{2}\in c(f_{1}f_{2}, t_{2}), t_{1}, t_{2}\in \mathbb{Z}^{3} \big\}.
\end{equation}

\end{exam}

\begin{prop} Each element of $V_{n}$ gives rise to a cluster strand.
\end{prop}
\begin{proof} For every element $v$ of $V_{n}$, one can find $r_{1},\ldots, r_{l}$  elements of $K(f_{1}, \ldots, f_{n})$    such that  $v$   can be written uniquely as follows
     \begin{equation}\label{}
        \nonumber v=r_{1}(f_{1}, \ldots, f_{n})y^{\beta_{11}}_{1,m_{11}}\cdots y^{\beta_{1n}}_{n,m_{1n}}+\ldots+r_{l}(f_{1}, \ldots, f _{n})y^{\beta_{l1}}_{1,m_{l1}}\cdots y^{\beta_{ln}}_{n,m_{ln}}.
     \end{equation}
      Such that a  $1-1$ map $\sigma: [1, l]\rightarrow \mathbb{Z}^{n}_{\geq0}\times \mathbb{Z}^{n}$ can be defined with $\sigma (j)=(\sigma _{1}(j), \sigma _{2}(j))$, where $\sigma _{1}(j)=(\beta _{j1},\ldots,\beta _{jn})$ and  $\sigma _{2}(j)=(m _{j1},\ldots,m _{jn}),j\in [1,l]$. Using Definition 5.8, one can introduce a cluster strand  $S_{l}(\sigma, r)$ with $r=(r_{1}, \ldots, r_{l})$ and $\sigma$ as defined above.
\end{proof}
 We denote the cluster strand associated to $v$ by  $S_{l}(\sigma, r)(v)$.\\
 \begin{que}
 Does the cluster strand $S_{l}(\sigma, r)(v)$ depend on the choices of $r$ or $\sigma$?
 \end{que}

The following  lemma and Remarks 5.14 provide some basic properties of the cluster strands.
\begin{lem}
\begin{enumerate}

\item Let  $v$ be an  element of $S_{l}(r, \sigma )$. Then  $S_{l}(r, \sigma)(v)=S_{l}(r, \sigma)$;

\item We have $\mathrm{S}_{l}(\sigma, g)=\mathrm{S}_{l}(\sigma, f)$  if and only if   for every $i \in [1, n]$, either  $g_{i} \in c(f_{i}, t_{i} )$ for some $t_{i} \in \mathbb{Z}^{n}$   or $f_{i}  \in c(g_{i}, t'_{i})$ for some $t'_{i} \in \mathbb{Z}^{n}$;
\item
Let $g=(g_{1},\ldots,g_{n})$ such that   $g_{i} \in c(f_{i}, t_{i})$  for some  $t_{i}\in \mathbb{Z}^{n},  i\in [1, n]$. Then $\mathrm{S}_{l}(\sigma', g)=\mathrm{S}_{l}(\sigma, f)$ if and only if $\forall j \in [1,l], \sigma_{1} '(j)=\sigma_{1} (j)$ and $\sigma_{2} '(j)=\sigma_{2} (j)+ q_{j} $ for some $q_{j}\in \mathbb{Z}^{n}$.
\end{enumerate}
\end{lem}
 \begin{proof}

For Part (1). Fix $v=\sum^{l}_{j=1}g_{j}y^{\beta_{j1}}_{1,m'_{j1}}\cdots y^{\beta_{jn}}_{n,m'_{jn}}\in S_{l}(\sigma, r )$. Then we must have, for every  $j\in [1, l], g_{j}=\alpha'_{j1}\cdots\alpha'_{jn}  r_{j}(\theta^{t' _{j1}} _{1}(f_{1}),\ldots,\theta^{t' _{jn}} _{n}(f_{n}))\in c(r_{j}, t'_{j})$, for some $t'_{j}=(t'_{j1},\ldots, t'_{jn})\in \mathbb{Z}^{n}$. A typical element of $c(r_{j}, t_{j})$ would be of the form $\alpha_{j1}\cdots\alpha_{jn}  r_{j}(\theta^{t_{j1}} _{1}(f_{1}),\ldots,\theta^{t _{jn}} _{n}(f_{n}))$ with $\alpha _{ji}\in M(f_{i}), \forall i\in [1, n]$, which can be written as
 \begin{equation}\label{}
  \nonumber  \alpha_{j1}\cdots\alpha_{jn} \alpha'_{jn}(\theta^{t_{j1}-t'_{j1}} _{1}(f^{-1}_{1}))\cdots \alpha'_{j1}(\theta^{t _{jn}-t' _{jn}} _{n}(f^{-1}_{n})) g_{j}(\theta^{t_{j1}-t'_{j1}} _{1}(f_{1}),\ldots,\theta^{t _{jn}-t' _{jn}} _{n}(f_{n}))
 \end{equation}
 which is an element of  $c(g_{j}, t_{j}-t'_{j})$.   Thus, any element  of the  following form $\sum^{l}_{j=1}r_{j}y^{\beta_{j1}}_{1,m_{j1}+t_{j1}}\cdots y^{\beta_{jn}}_{n,m_{jn}+t_{jn}}$ is in fact an element of $S_{l}(g, \sigma')$, where $\sigma_{1}'(j)=\sigma_{1}(j)$ and $\sigma_{2}'(j)=\sigma_{2}(j)+t_{j}, j \in [1, l]$. Then $S_{l}(\sigma, r)\subseteq S_{l}(\sigma ', g)$. But from the Proof of Proposition 5.10, one can see that  $S_{l}(\sigma', g)=S_{l}(\sigma', g)(v)$.  Again from the proof of Proposition 5.10, one can see that $S_{l}(\sigma, g)(v)\subseteq S_{l}(\sigma, r)$. Therefore,  $S_{l}(\sigma, r )=S_{l}(\sigma, g)(v)$.

For Part (2).  $(\Rightarrow)$ is Obvious. For the other direction $(\Leftarrow)$. Without loss of  generality, let $g_{j}\in  c(f_{j}, t_{j})$ for some $t_{j}=(t_{j1}, \ldots, t_{jn})\in \mathbb{Z}^{n}$. Then for every $j\in [1, l]$, there are $\alpha_{ji}\in M(f_{i}), i\in [1, n]$ such that $g_{j}=\alpha_{j1}\cdots\alpha_{jn}  f_{j}(\theta^{t _{j1}} _{1}(f_{1}),\ldots,\theta^{t _{jn}} _{n}(f_{n}))$. Now, let $v\in S_{l}(g, \sigma)$. Hence, we have
\begin{eqnarray*}
% \nonumber to remove numbering (before each equation)
 \nonumber v &=& \sum^{l}_{j=1}g_{j}y^{\beta_{j1}}_{1,m'_{j1}}\cdots y^{\beta_{jn}}_{n,m'_{jn}} \\
 \nonumber  &=& \sum^{l}_{j=1}\alpha_{j1}\cdots\alpha_{jn}  f_{j}(\theta^{t _{j1}} _{1}(f_{1}),\ldots,\theta^{t _{jn}} _{n}(f_{n}))y^{\beta_{j1}}_{1,m'_{j1}}\cdots y^{\beta_{jn}}_{n,m'_{jn}}\in S_{l}(f, \sigma).
\end{eqnarray*}
Therefore $S_{l}(\sigma, f)=S_{l}(\sigma, g)(v)=S_{l}(\sigma, g)$ thanks to Part (1) of this lemma.

 For Part (3). First for $(\Rightarrow)$. One can  see that, if $\sigma '_{1}(j)=\sigma _{1}(j), \forall j\in [1,l]$, then $\sigma' = \sigma+(0,q), q\in \mathbb{Z}^{n}$. Now, assume that  $\sigma '(j_{0})\neq \sigma(j_{0})+(0,q_{j})$ for some $j_{0}\in [1,l]$ and for every $q\in \mathbb{Z}^{n}$. Hence $\sigma _{1}'(j_{0}) \neq \sigma _{1}(j_{0})$. Then the element
 \begin{equation}\label{}
  \nonumber v_{0}= g_{j_{0}}y^{\beta_{j_{0}1}}_{1,m_{j_{0}1}+t_{j_{0}1}}\cdots y^{\beta_{j_{0}n}}_{n,m_{j_{0}n}+t_{j_{0}n}}+\sum^{l}_{j\in [1, l]\setminus \{j_{0}\}}g_{j}y^{\beta_{j1}}_{1,m_{j1}+t_{j1}}\cdots y^{\beta_{jn}}_{n,m_{jn}+t_{jn}}
 \end{equation}
 is an element of $\mathrm{S}_{l}(\sigma', g)$  with $\sigma' _{1}(j)=(\beta_{j1},\ldots,\beta_{jn})$. However, $v_{0}$  is not an element of $\mathrm{S}_{l}(\sigma, f)$. $(\Leftarrow)$ is immediate.

\end{proof}

\begin{defn} Any submodule of $V_{n}$ generated by a cluster strand  $\mathrm{S}_{l}(\sigma, r)$ is called a strand  submodule and is denoted by $W_{l}(\sigma, r)$.\end{defn} In the occasions, when we want to emphasis  on a certain element $v$ of $V_{n}$, we will denote the  strand submodule associated to the cluster strand $\mathrm{S}_{l}(\sigma, r)(v)$ by $W_{l}(\sigma, r)(v)$ or, for the sake of simplicity, by $W_{l}(w)$.

Let $\mathcal{M}(E)$ be the set of all monomials formed from  elements of  the set $E=\{x_{1},\ldots,x_{n},y_{1},\ldots,y_{n}\}$.  \emph{A special  cluster strand} is defined to be a subset of a  full cluster strand $S_{l}(\sigma, r)$ of the form
\begin{equation}\label{}
  \widehat{S}_{l}(\sigma, r):=\Big\{\sum^{l}_{j=1}g_{j}y^{\beta_{j1}}_{1,m_{j1}+t_{1}}\cdots y^{\beta_{jn}}_{n,m_{jn}+t_{n}}| t=(t_{1},\ldots, t_{n})\in \mathbb{Z}^{n}, g_{j}\in c(r_{j}, t), j \in [1, l]\Big\}.
\end{equation}

The submodule of $W_{l}(\sigma, r)$ generated by the special cluster strand $\widehat{S}_{l}(\sigma, r)$ is called  \emph{special stand module} and is denoted by $\widehat{W}_{l}(\sigma, r)$.

\begin{rems}\begin{enumerate}
\item
 Let $\widehat{S}_{l}(\sigma, h)$ be a special cluster strand. Then
 \begin{enumerate}
 \item
\begin{equation}\label{}
  \nonumber \widehat{S}_{l}(\sigma , h)  \ \text{is a proper subset of} \ S_{l}(\sigma , h);
\end{equation}
\item
\begin{equation}\label{}
  \nonumber \widehat{S}_{l}(\sigma , h)=\mathcal{M}(E)w, \ \text{for every } \ w \in \widehat{S}_{l}(\sigma, h).
  \end{equation}
 \end{enumerate}
\item There is a  bijection between the set of all  cyclic submodules of $V_{n}$ and the set of all special strand submodules.
\end{enumerate}
\end{rems}
\begin{proof}Part (1) is straight forward. For Part (2), let $W$ be a cyclic module generated by $w$ with associated cluster strand $S_{l}(\sigma, r)(w)$. Then by the definition of special cluster strands, we have $\widehat{W}_{l}(\sigma, r)$ is a submodule of $W$. One can realize that $W$ is a submodule of $\widehat{W}_{l}(\sigma, r)$ too,  if we recall that $W$ is cyclic module generated by $w$ which is an element of $\widehat{S}_{l}(\sigma, r)$. The bijection is defined to send  $W$  to $\widehat{S}_{l}(\sigma, r)$.
\end{proof}

\begin{prop}
\begin{enumerate}
\item
Every strand submodule $W_{l}(f, \sigma)$ can be identified with a  sum of (identical) copies of the cluster strand $S_{l}(h, \sigma)$.
\item Every submodule $W$ of $V_{n}$ is a sum of some strand submodules. In particular, $W$ is generated by a set of cluster
    strands.
    \end{enumerate}
\end{prop}
\begin{proof}
\begin{enumerate}
\item
First we show that the extensions of the action of the elements of $\mathcal{M}(E)$, induced by  (5.3) and (5.4),  keeps the cluster strands invariant. One can see that for every  $g\in K(f_{1}, \ldots, f_{n})$ and $t\in \mathbb{Z}^{n}$ the coefficients set $c(g, t)$ is invariant under the actions of the elements of $E$ and then under elements of $\mathcal{M}(E)$. Now, let $v=\sum^{l}_{j=1}g_{j}y^{\beta_{j1}}_{1,m_{j1}}\cdots y^{\beta_{jn}}_{n,m_{jn}}$ be an element of the cluster strand $S_{l}(h, \sigma)$. Recalling that the actions given in (5.3) and (5.4) define a fully faithful representation, one can see that under the actions of elements of $E$ the length $l$ stays unchanged with respect to $h$ and  $\sigma$, which will stay unchanged too. Hence for any monomial $m\in \mathcal{M}(E)$, we have $m(v)\in S_{l}(h, \sigma)$, more precisely
    \begin{equation}\label{}
       m(v)=\sum^{l}_{j=1}\alpha_{j1}\cdots \alpha_{jn}g_{j}(\theta^{t_{j1}}_{1}(f_{1}),\ldots, \theta^{t_{jn}}_{1}(f_{n}))y^{\beta_{j1}}_{1,m_{j1}+t_{j1}}\cdots y^{\beta_{jn}}_{n,m_{jn}+t_{jn}}.
    \end{equation}
    Where $\alpha_{ji}\in M(f_{i})$ and $t_{ji}\in \mathbb{Z},$ $\forall i \in [1, n], j\in [1, l]$. Recall that, elements of  $W_{l}(h, \sigma)$ are finite sums of finite products of elements of $R(\theta, \xi, n)$ acting on an element of  $S_{l}(f, \sigma)$. But every element of $R(\theta, \xi, n)$ can be written as a $K(f_{1},\ldots, f_{n})$-linear combination of elements of $\mathcal{M}(E)$. Then from (5.9) elements of $W_{l}(h, \sigma)$ are finite sum of elements of $S_{l}(h, \sigma)$. In the same time one can obviously see that every sum of elements of $S_{l}(h, \sigma)$ must be an element of $W_{l}(h, \sigma)$.

\item
 We first notice that, from Part (1) of Lemma 5.12, we conclude that every two cluster strands are either identical or have zero intersection. So we can
     introduce the following equivalence relation on  $V_{n}$
\begin{equation}\label{}
   \forall s, s'\in V_{n}, \  s \sim s'  \ \text{if and only if} \  s  \ \text {and} \  s' \  \text {belong to the same cluster strand}.
\end{equation}
Let $W$ be a submodule of $V_{n}$ and $W^{*}=W/\sim $. Here, every $w^{*}\in W^{*}$, is the intersection of $W$ with the cluster strand $S_{l}(\sigma, f)(w)$. If $W^{*}_{l}(w)$ denote the submodule of $W$  generated by $w^{*}$. Then we have the
following identity
\begin{equation}\label{}
    W=\sum^{} _{w^*\in W^*} W^{*}_{l}(w).
\end{equation}

    \end{enumerate}
\end{proof}

\begin{prop}
Let $S_{l}(\sigma, h)$ be a full cluster strand. Then any two  strand submodules of  $W_{l}(\sigma, h)$  have a non-zero intersection.

 \end{prop}

\begin{proof}
  Let $W_{1}=W_{l_{1}}(\sigma^{\mathfrak{1}}, h^{\mathfrak{1}})(w_{1})$ and $W_{2}=W_{l_{2}}(\sigma^{\mathfrak{2}}, h^{\mathfrak{2}})(w_{2})$  be any two proper strand submodules of $W_{l}(\sigma, h)$. From Proposition 5.15 and Proposition 5.10, one can see that the cluster stands $S_{l_{1}}(\sigma^{\mathfrak{1}}, h^{\mathfrak{1}})(w_{1})$ and  $S_{l_{2}}(\sigma^{\mathfrak{2}}, h^{\mathfrak{2}})(w_{2})$  satisfy the following

     \begin{itemize}
       \item  There are two natural numbers $d_{1}$ and $d_{2}$ such that $l_{i}=d_{i}l, i= \mathfrak{1, 2}$;
       \item $ h^{\mathfrak{1}}=(h_{11},\ldots, h_{1l}, \ldots, h_{d_{\mathfrak{1}}1},\ldots, h_{d_{\mathfrak{1}}l})$ and $h^{\mathfrak{2}}=(h'_{11},\ldots, h'_{1l}, \ldots, h'_{d_{\mathfrak{2}}1},\ldots, h'_{d_{\mathfrak{2}}l})$ where there are  $t_{j_{2}}\in \mathbb{Z}^{n}$ such that $h_{j_{1}j_{2}}, h'_{j_{1}j_{2}} \in c(h_{j_{2}}, t_{j_{2}})$ for every $j_{1}\in [1, d_{i}],  j_{2}\in [1, l]$;
       \item For $i=\mathfrak{1}, \mathfrak{2}$ we have  $\sigma^{i}:\{11,12, \ldots, 1l, \ldots, d_{i}1, \ldots, d_{i}l \}\rightarrow \mathbb{Z}_{>0}^{n}\times \mathbb{Z}^{n}$, such that $\sigma^{i}=(\sigma^{i}_{1}, \sigma^{i}_{2})$ where $\sigma^{i}_{1}(j_{1}j_{2})=\sigma_{1}(j_{2})$ and $\sigma^{i}_{2}(j_{1}j_{2})=\sigma_{2}(j_{2})+t_{j_{1}j_{2}}$  for some $t_{j_{1}j_{2}}\in \mathbb{Z}^{n}$, for all $j_{2}\in [1, l]$.
     \end{itemize}
     Now we will show that the sum of any $d_{i}$-elements of $S_{l}(\sigma, h)$ is an element of $S_{l_{i}}(\sigma^{i}, h^{i})(w_{i})$, for $i=\mathfrak{1, 2}$. Consider the following two elements
     \begin{equation}\label{}
     \nonumber  w_{j_{1}}=\sum_{j_{2}=1}^{l}h_{j_{1}j_{2}}y^{\beta_{j_{2}1}}_{1,m_{j_{2}1}+t^{1}_{j_{1}j_{2}}}\cdots y^{\beta_{j_{2}n}}_{n,m_{j_{2}n}+t^{n}_{j_{1}j_{2}}}
     \end{equation}

    and
\begin{equation}\label{}
     \nonumber w'_{j_{1}}=\sum_{j_{2}=1}^{l}h'_{j_{1}j_{2}}y^{\beta_{j_{2}1}}_{1,m_{j_{2}1}+t^{1}_{j_{1}j_{2}}}\cdots y^{\beta_{j_{2}n}}_{n,m_{j_{2}n}+t^{n}_{j_{1}j_{2}}}.
\end{equation}
One can see that the elements $w_{j_{1}}$ and  $w'_{j_{1}}$ are elements of $S_{l}(\sigma, h)$ for every $j_{1}\in [1, d_{i}]$ for $i=1, 2$. Then the cluster strands associated to $w_{j_{1}}$ and  $w'_{j_{1}}$ coincide with $S_{l}(\sigma, h)$ for every $j_{1}\in [1, d_{i}], i=1,2$, thanks to Part (1) of Lemma 5.12.

for every $j_{1}\in [1, d_{i}]$, we have $s\in S(w_{j_{1}})=S(w'_{j_{1}})$ for every $s\in S(w)$.  Let $l'$ be the least common multiple of $l_{1}$ and $l_{2}$. So, $l'=n_{i}l_{i}$,
     for some $n_{i}\in \mathbb{N}, i=1,2$. Consider the element
\begin{equation}\label{}
  \nonumber  w'=\sum_{i=1}^{l'}s_{i}, \text {where} \ \ s_{i} \in S(w) \setminus \{s_{1}, \ldots, s_{i-1}\}, \forall i\in [1, l'].
\end{equation}
One can see that, $w'$ is in deed a sum of $n_{i}d_{i}$-elements of the cluster strand of $w_{j_{i}}$  and  elements of $S_{l_{i}}(\sigma^{i}, h^{i})(w_{i})$ are sums of $d_{i}$-elements of the cluster strand  $w_{j_{i}}, i=1,2$. Then $w'$ is a sum of $d_{i}$-elements of the cluster stared $S_{l_{i}}(\sigma^{i}, h^{i})(w_{i}), i=1, 2$.  Therefore from  Part (1) of Proposition 5.15, we have
\begin{equation}\label{}
\nonumber w' \in W(w_{1})\cap W(w_{2}).
\end{equation}
 \end{proof}
 The following corollary is  a consequence of the proof of Proposition 5.16. Let $S_{l}(\sigma, h)$ be a cluster strand with a strand module $W_{l}(\sigma, h)$. For every natural number  $j$ we introduce a subset of $S_{l}(\sigma, h)$ give by
  \begin{equation}\label{}
    \nonumber S^{j}_{l}(\sigma, h)=\{s_{1}+\cdots+s_{j};  \ \ s_{i} \in S_{l}(\sigma, h)\setminus\{s_{1}, \ldots, s_{i-1}\}, \forall i\in [1, j]\}.
  \end{equation}

 \begin{cor}
  \begin{enumerate}
 \item For every  $w^{j}, s^{j}\in S^{j}_{l}(\sigma, h)$, the cluster strands $S_{jl}(\sigma^{\mathfrak{1}}, h^{\mathfrak{2}})(s^{j})$ and $S_{jl}(\sigma^{\mathfrak{2}}, h^{\mathfrak{2}})(w^{j})$, defined in the Proof of Proposition 5.16 are coincide  and with length of $jl$.

 \item Let $W(j)$ denote the strand module of $S_{jl}(\sigma^{\mathfrak{1}}, h^{\mathfrak{2}})(s^{j})$. Then we have the following descending chain of strand modules
 \begin{equation}\label{}
  \nonumber  W_{l}(\sigma, h)\supseteq W(j)\supset W(2j)\supset \ldots \supset W(nj)\supset \ldots, \forall j \in \mathbb{N}.
 \end{equation}

 \end{enumerate}
 \end{cor}
 % -----------------------------------------------------------------------------

\subsection*{Acknowledgment}
I am grateful to an anonymous referee for the very useful suggestions and comments that helped a lot in getting this paper in its final form. This work has  been started while I was attending the topics in Mathematics class taught by late Alexander L. Rosenberg. Many thanks to Zongzhu Lin for all valued discussions.
% -----------------------------------------------------------------------------

\end{document}